\documentclass[a4paper, english,reqno, 11pt]{amsart}
\usepackage{amsmath,amsthm,amssymb}
\usepackage{enumitem}
\usepackage{float}

\usepackage{ifxetex,ifluatex}
\if\ifxetex T\else\ifluatex T\else F\fi\fi T%
  \usepackage{fontspec}
\else
  \usepackage[T1]{fontenc}
  \usepackage[utf8]{inputenc}
  \usepackage{lmodern}
\fi

\usepackage[margin=1.25in]{geometry}

\usepackage{tikz}
\usetikzlibrary{shapes.geometric}
\usetikzlibrary{decorations.pathreplacing}
\usepackage{xcolor}
\usepackage{hyperref}
\usepackage{cleveref}

\usepackage{aliascnt}
\newcommand{\newsharedtheorem}[2]{%
    \newaliascnt{#1}{theorem}%
    \newtheorem{#1}[#1]{#2}%
    \aliascntresetthe{#1}%
}

\theoremstyle{plain}
\newtheorem{theorem}{Theorem}[section]
\newsharedtheorem{corollary}{Corollary}
\newsharedtheorem{lemma}{Lemma}
\newsharedtheorem{proposition}{Proposition}

\theoremstyle{definition}
\newsharedtheorem{conjecture}{Conjecture}
\newsharedtheorem{example}{Example}
\newsharedtheorem{remark}{Remark}

\newcommand{\zero}{\hat 0}
\newcommand{\one}{\hat 1}
\DeclareMathOperator{\Lab}{Lab}
\DeclareMathOperator{\des}{des}
\DeclareMathOperator{\Des}{Des}
\DeclareMathOperator{\rk}{rk}
\newcommand{\closure}{\operatorname{cl}}
\newcommand{\Chow}{\operatorname{H}}
\newcommand{\lowvalley}{\operatorname{lv}}

\newcommand{\aaa}{\mathbf{a}}
\newcommand{\bbb}{\mathbf{b}}
\newcommand{\ZZ}{\mathbb{Z}}
\newcommand{\ZZab}{\ZZ\langle \aaa,\bbb \rangle}
\newcommand{\NN}{\mathbb{N}}
\newcommand{\RR}{\mathbb{R}}

\newcommand{\cF}{\mathcal{F}}
\newcommand{\cL}{\mathcal{L}}

\newcommand{\cP}{\mathcal{P}}

\newcommand{\cS}{\mathcal{S}}

\usepackage{microtype}
\newcommand{\lp}[1]{\text{\scalebox{0.85}{\textls[-60]{#1}}}}
\newcommand{\LPM}[2]{M[\lp{#1},\lp{#2}]}

\newcommand{\Dfn}[1]{\emph{\bfseries #1}}

\NewDocumentCommand{\drawcells}{O{gray!50}m}{
    \foreach \x/\y in {#2}
        \draw[#1] (\x, \y) rectangle ++(1,1);
}

\NewDocumentCommand{\drawpath}{sO{thick}m}{
    \coordinate (p) at (0,0);
    \foreach \s [count=\k] in {#3} {
        \if\s N \draw[#2] (p) -- ++(0,1) \IfBooleanT{#1} {node[midway,fill=white,inner sep=1pt,font=\scriptsize]{$\k$}} coordinate (p);
        \else \draw[#2] (p) -- ++(1,0) \IfBooleanT{#1} {node[midway,fill=white,inner sep=1pt,font=\scriptsize]{$\k$}} coordinate (p);
        \fi
    }
}

\colorlet{gsT}{green!50!black!35}
\colorlet{gsC}{orange!40}
\colorlet{gsP}{violet!35}
\tikzset{
    graph label/.style={fill=white, inner sep=1pt, font=\scriptsize},
    e/.style={},
    lv/.style={graph label/.append style={ellipse, draw=red, thick, inner sep=0.4pt}},
    T/.code={\graphmarker{gsT}},
    C/.code={\graphmarker{gsC}},
    P/.code={\graphmarker{gsP}},
    lvT/.style={lv,T},
    dots/.style={dotted, thick},
    highlighter gsP/.style={dash pattern=on 7pt off 3.5pt, line cap=butt},
}
\newcommand{\graphmarker}[1]{\def\graphmarkerTmp{#1}\ifx\graphmarkerTmp\graphmarkerPass\gdef\graphmarkerHit{1}\fi}

\NewDocumentCommand{\drawgraph}{O{}mm}{
    \coordinate (v0) at (0,0);
    \foreach \i in {1,...,#2}
        \coordinate (v\i) at ({1.2*(\i-(#2+1)/2},{1.5});
    \foreach \e/\a/\b/\bend/\style in {#3} {
        \tikzset{graph edge/.style/.expand once=\style}
        \ifx\e\empty
            \draw[graph edge] (v\a) to[bend left=\bend] (v\b);
        \else
            \draw[graph edge] (v\a) to[bend left=\bend] node[graph label] {$\e$} (v\b);
        \fi
    }
    \foreach \i in {1,...,#2}
        \fill (v\i) circle (1.5pt);
    \filldraw[fill=white] (v0) circle (2pt);
    \if\relax\detokenize{#1}\relax\else
        \foreach \i/\t in {#1} {
            \ifnum\i=0 \node[below=1pt, font=\scriptsize] at (v0) {$\t$};
            \else \node[above=1pt, font=\scriptsize] at (v\i) {$\t$};
            \fi
        }
    \fi
    \foreach \m in {gsT,gsC,gsP} {
        \edef\graphmarkerPass{\m}
        \begin{scope}[blend mode=multiply]
        \begin{scope}[transparency group]
            \foreach \e/\a/\b/\bend/\style in {#3} {
                \gdef\graphmarkerHit{0}
                {\tikzset{graph edge/.style/.expand once=\style, graph edge}}
                \ifnum\graphmarkerHit=1
                    \draw[\graphmarkerPass, line width=5.5pt, line cap=round,
                        highlighter \graphmarkerPass/.try]
                        (v\a) to[bend left=\bend] (v\b);
                \fi
            }
        \end{scope}
        \end{scope}
    }
    \def\graphmarkerPass{}
}

\tikzset{snakestyle/.style={fill=blue!15, draw=gray!50}}
\newcommand{\drawM}{
\drawcells{0/0, 1/0, 0/1, 1/1, 2/1, 1/2, 2/2}
\drawpath{N,N,E,N,E,E}
\drawpath{E,E,N,E,N,N}
}

\title{\texorpdfstring{$\aaa\bbb$}{ab}-indices for lattice path matroids}

\author[C.~Aydin]{Cenk Aydin}
\author[G.~Caiolo]{Giovanni Caiolo}
\author[C.~Giardino]{Chiara Giardino}
\author[M.~Gruß]{Mara Gruß}
\author[M.~H\"aglund]{Mikael H\"aglund}
\author[E.~Hoster]{Elena Hoster}
\author[C.~Stump]{Christian Stump}
\address[C.~Aydin \& C.~Giardino \& C.~Stump]{Ruhr-Universit\"at Bochum, Germany}
\email{\{cenk.aydin,chiara.giardino,christian.stump\}@rub.de}
\address[G.~Caiolo]{Scuola Normale Superiore di Pisa, Pisa, Italy}
\email{giovanni.caiolo@sns.it}
\address[M.~Gruß]{Freie Universität Berlin, Germany}
\email{mara.gruss@fu-berlin.de}
\address[M.~H\"aglund]{KTH Royal Institute of Technology, Sweden}
\email{mihaglu@kth.se}
\address[E.~Hoster]{Università di Bologna, Italy}
\email{elena.hoster@unibo.it}

\date{\today}
\usepackage[backend=bibtex,bibencoding=utf8,style=alphabetic,maxnames=5]{biblatex}
\begin{document}
\begin{abstract}
    The base polytopes of connected lattice path matroids admit subdivisions into \emph{snakes}.
    We use these to obtain an explicit positive combinatorial description of their $\aaa\bbb$-indices, their Chow polynomials and their chain polynomials.
    We finally apply this description to the particular case of the $\aaa\bbb$-index of Schubert matroids.
\end{abstract}

\maketitle
\section{Introduction and main results}
\subsection{Motivation}\label{subsec:intro-motivation}
Let~$M$ be a matroid on the ground set~$[n]=\{1<\dots< n\}$ of rank~$r$, and let $\cL = \cL(M)$ be its lattice of flats.
For a maximal chain $\cF=\{F_0 \prec F_1 \prec \dots \prec F_r\}$ of flats in~$\cL$, define the sequence $\lambda(\cF)=(\lambda_1,\dots , \lambda_r)$ with  $\lambda_i = \min( F_{i}\setminus F_{i-1} )$ and denote the set of \Dfn{$M$-sequences} arising this way by
\[
  \Lab(M) = \big\{ \lambda(\cF) \mid \cF \text{ maximal chain } \big\}.
\]
Such sequences are \emph{minimal labelings} of geometric lattices
and closely related to the theory of broken circuits and lexicographic shellings; see, for example, \cite{Rota1964,Bjorner1980,Bjorner1992}.
In particular, decreasing label sequences recover the classical \emph{NBC bases}.
Observe that not every ordering of an NBC basis is an $M$-sequence, and $M$-sequences are called NBC\(^{+}\)-bases in~\cite{HershSundaram2026}.

\smallskip
Towards our main results, we start with a characterization of such NBC\(^{+}\)-bases.
For an integer sequence $\lambda = (\lambda_1 , \dots , \lambda_k )$ and a subset $C\subset \ZZ$, write $\lambda|_C$ for the subsequence of entries of $\lambda$ lying in $C$.

\begin{proposition}\label{prop:LabM-general}
    Let $\lambda=(\lambda_1,\dots , \lambda_r)$ be an integer sequence.
    Then $\lambda$ is an $M$-sequence, $\lambda\in \Lab(M)$, if and only if the following two properties hold:
    \begin{itemize}
        \item \textbf{Ordered basis:}
        The set $B_\lambda = \{\lambda_1,\dots , \lambda_r\}$ is a basis of~$M$.
        \item \textbf{Fundamental circuits:}
        For every non-loop $x \in [n] \setminus B_\lambda$, the unique circuit~$C_x \subseteq B_\lambda \cup\{x\}$ of $M$ satisfies that~$\lambda|_{C_x}$ ends with an integer smaller than~$x$.
    \end{itemize}
\end{proposition}

One motivation for studying $M$-sequences comes from the following relation to flag $h$-vectors.
Following \cite{Stanley1974}, the \Dfn{flag $h$-vector}~$\beta_\cL$ of~$\cL$ in a subset~$T\subseteq [r-1]$ is
\begin{equation*}
    \beta_\cL(T) = \# \{ \lambda\in \Lab(M) \mid \Des(\lambda) = T \} \,,
\end{equation*}
where~$\Des(\lambda) = \{ i \mid \lambda_i > \lambda_{i+1} \}$ is the set of \Dfn{descents} of $\lambda$.
We encode it in the \Dfn{$\aaa\bbb$-index} $\Psi_\cL(\aaa,\bbb) \in \ZZ\langle \aaa,\bbb\rangle$, a polynomial in two non-commuting variables $\aaa,\bbb$ defined by
\begin{equation}\label{eq:def-ab-index}
    \Psi_\cL(\aaa,\bbb)
    = \
     \sum_{T} \beta_\cL(T) \cdot \operatorname{m}_T(\aaa,\bbb) \,,
     \quad 
     \operatorname{m}_T(\aaa,\bbb) = m_1 \cdots m_{r-1}\,,
     \quad
     m_i = 
     \begin{cases}
         \bbb \,,   & i\in T \,, \\
         \aaa \,,   & i\notin T \,,
     \end{cases}
\end{equation}
where the sum ranges over all subsets~$T\subseteq\{1,\dots , r-1\}$.
We may also call $\Psi_\cL(\aaa,\bbb)$ the $\aaa\bbb$-index of $M$ and write $\Psi_M(\aaa,\bbb) = \Psi_\cL(\aaa,\bbb)$ and $\beta_M = \beta_\cL$.

\smallskip
Based on \Cref{prop:LabM-general}, the main result of this paper is \textit{an intrinsically positive combinatorial description} of the $\aaa\bbb$-index of lattice path matroids.
This simultaneously yields combinatorial formulas for all invariants that are determined by the flag $h$-vector~$\beta_\cL(T)$.
These include, for example, the flag~$f$-vector of $\cL$, the extended $\aaa\bbb$-index of $\cL$, the $f$-/$h$-vector of the order complex of $\cL$, the Chow polynomial of $M$, and many more.

\subsection{Background on lattice path matroids}\label{subsec:intro-background-lpm}
A lattice path $\gamma=\gamma_1\dots \gamma_n$ is a path in~$\RR^2$ along the integer grid, going in north steps $\lp{N}=+(0,1)$ and east steps $\lp{E}=+(1,0)$.
Let~$U$ and~$L$ be two lattice paths, both going from $(0,0)$ to $(n-r,r)$ such that $U$ stays weakly above $L$.
The \Dfn{lattice path matroid} (\Dfn{LPM}) $M=M[U,L]$ is the matroid on ground set $[n]$ and of rank $r$ whose bases are determined by lattice paths weakly between the \Dfn{upper path} $U$ and the \Dfn{lower path} $L$. That is, a set $B\subseteq[n]$ is a basis of $M$ if the path
\[
    P=p_1\cdots p_n\,, \text{with }
    p_i=
    \begin{cases}
        \lp{N} & \text{if } i\in B \,, \\
        \lp{E} & \text{if } i\notin B \,,
    \end{cases}
\]
never goes above~$U$ and never goes below~$L$.
A common east step (resp. north step) of~$U$ and~$L$, that is, a step taken by both paths from the same lattice point, corresponds to a loop (resp. coloop) of~$M$.
If $U$ and $L$ have, beside the endpoints, no lattice points in common, we call $M[U,L]$ \Dfn{connected}.
(In particular, loops and coloops are connected.)
Every lattice path matroid $M=M[U,L]$ can be decomposed into a direct sum $M= M_1 \oplus \cdots \oplus M_m$ of connected LPMs. 
Ordered in the obvious way, this decomposition is unique.
In other words, the upper and lower path of $M$ can be uniquely written as
\[
    \begin{array}{rcl@{\qquad}l}
    U = U_1\cdots U_m\,, & \\[-1ex]
                        & \text{such that } M_i = M[U_i,L_i] \text{ are connected for all } i\,. \\[-1ex]
    L = L_1\cdots L_m\,, &
    \end{array}
\]
In this case, $M_i$ are the connected components of $M$ and $c(M)=m$ denotes the number of components.

The region between $U$ and $L$ is called the \Dfn{diagram} of $M$.
We say that an LPM~$N$ on~$[n]$ is \Dfn{contained} in~$M$ if every basis of $N$ is also a basis of $M$.
Equivalently, $N$~is contained in~$M$ if they have the same rank and ground set and the diagram of~$N$ is contained in the diagram of~$M$.
Finally, a \Dfn{snake matroid} (or \Dfn{snake}) $S$ is an LPM which has neither a loop nor a coloop and for which the diagram of $S$ does not contain a $(2\times 2)$-box\footnote{In the literature, snakes are usually required to be connected. We do not assume this additional property to necessarily hold.}.

\begin{figure}[H]
    \centering
    \begin{tikzpicture}[scale=0.75]
        \drawcells[snakestyle]{0/0, 1/1, 1/2, 2/2}
        \drawM
    \end{tikzpicture}
    \caption{An LPM containing a (non-connected) snake, shaded in blue.}
    \label{fig:running-M}
\end{figure}
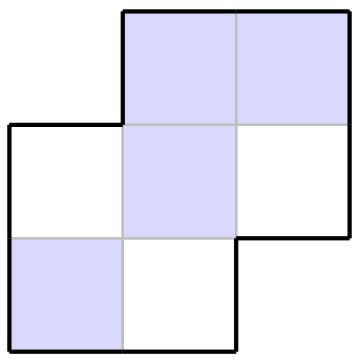

\begin{example}\label{ex:running-M}
    Throughout this paper, we use the LPM $M=\LPM{NNENEE}{EENENN}$ as a running example.
    It contains the snake~$S=\LPM{NENNEE}{ENENEN}$ shaded in \Cref{fig:running-M}.
\end{example}

\subsection{Main results on lattice path matroids}\label{subsec:intro-main-results-lpm}
Bidkhori \cite{Bidkhori2010} showed that the base polytope $\mathcal P(M)$ of a connected lattice path matroid $M$ carries a canonical subdivision into snakes contained in $M$.
This \emph{snake subdivision} has been used to study several enumerative and geometric properties of LPMs, see \cite{Bidkhori2010,FerroniMoralesPanova2026, ChatelainRamirezAlfonsin2011, BenedettiKnauerValenciaPorras2023}.

Since the $\aaa\bbb$-index is a \emph{valuative invariant}, see \cite[Theorem 8.2]{ferroni2024valuative}, it can be computed from the snake subdivision of the base polytope $\mathcal P(M)$ as follows.
Let~$M=M[U,L]$ be a connected LPM.
Then
\begin{equation}\label{eq:ab-LPMs}
    \Psi_M(\aaa,\bbb)
    = \sum_{Q} (-1)^{c(Q)-1} \cdot \Psi_Q(\aaa,\bbb) \,,
\end{equation}
where the sum ranges over all snakes~$Q$ contained in~$M$, with the additional property that~$M[U_Q,L]$ and~$M[U,L_Q]$ are both connected.
We call such snakes~$Q$ \Dfn{internal faces of~$M$}, since their base polytopes are the internal faces of the base polytope $\mathcal P(M)$ under the snake subdivision.
Observe that the sign in~\eqref{eq:ab-LPMs} depends on the number of components in the decomposition of the snake $Q = S_1 \oplus \dots \oplus S_{c(Q)}$ into connected snakes.

\begin{figure}[H]
    \centering
    \begin{tikzpicture}[scale=0.75]
        \drawcells[snakestyle]{0/0, 1/1, 2/1, 2/2}
        \drawM
    \end{tikzpicture}
    \qquad
    \begin{tikzpicture}[scale=0.75]
        \drawcells[fill=gray!25, draw=gray!50]{0/0, 0/1, 1/2, 2/2}
        \drawM
    \end{tikzpicture}
    \caption{An example (left) and a non-example (right) of an internal face of an LPM.}
    \label{fig:running-faces}
\end{figure}
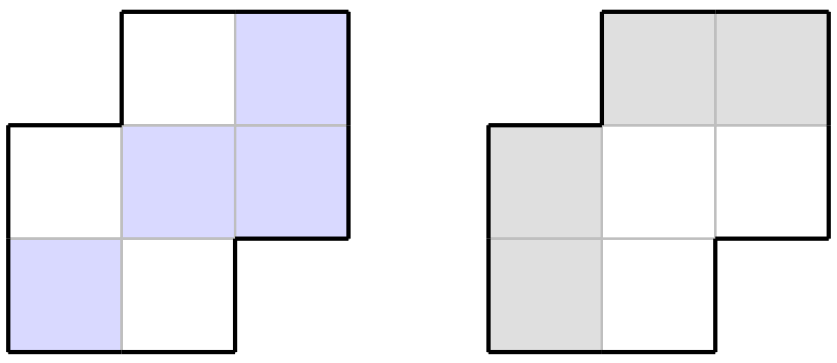

\begin{example}\label{ex:running-faces}
    Continuing \Cref{ex:running-M}, the left of \Cref{fig:running-faces} shows the internal face $Q = \LPM{NENENE}{ENEENN}$ of $M$ with its two components $\LPM{NE}{EN}$ and $\LPM{NENE}{EENN}$.
    The snake $Q' = \LPM{NNENEE}{ENNEEN} = \LPM{NNE}{ENN} \oplus \LPM{NEE}{EEN}$ on the right is \emph{not} an internal face because its upper path is $U$, so that $M[U,L_{Q'}] = Q'$ is not connected.
\end{example}
For a connected snake $S=M[U_S,L_S]$ contained in~$M$, we call an $S$-sequence~$\lambda$ \Dfn{$(M,S)$-surviving} if~$\lambda$ is \emph{not} a $Q$-sequence, $\lambda \notin \Lab(Q)$, for every non-connected internal face~$Q$ of~$M$ having the same upper path as~$S$, $U_Q = U_S$.
Note here and below that~$S$ is the unique connected snake having the same upper path as~$Q$.
The following theorem describes the $\aaa\bbb$-index as a sum over surviving sequences.

\begin{theorem}\label{thm:intro-ab-via-MS}
    The $\aaa\bbb$-index of a connected lattice path matroid $M$ equals
    \begin{equation*}
        \Psi_M(\aaa,\bbb)
        = \sum_{S } \sum_{\lambda } \operatorname{m}_{\Des(\lambda)}(\aaa,\bbb)\,,
    \end{equation*}
    where the sum ranges over all connected snakes~$S$ contained in~$M$ and all $(M,S)$-surviving sequences $\lambda\in \Lab(S)$.
\end{theorem}

Every positive specialization of the $\aaa\bbb$-index inherits a manifestly positive formula from this description.
We highlight two important consequences, which are both immediate from the previous theorem.

\begin{corollary}\label{cor:flag-h-LPM-via-snakes}
The flag $h$-vector of a connected lattice path matroid $M$ of rank $r$ is
    \begin{equation*}
        \beta_M(T)
        = \#\biggl\{
            (S,\lambda) \ \biggm| \
            \substack{ S \text{ connected snake contained in $M$, } \\
                \lambda \text{ is $(M,S)$-surviving, and }
                \Des(\lambda) = T} \ 
            \biggr\}  \,.
    \end{equation*}
\end{corollary}
\begin{proof}
    This follows from \Cref{thm:intro-ab-via-MS} by comparing coefficients with the monomial expansion of the $\aaa\bbb$-index in~\eqref{eq:def-ab-index}.
\end{proof}

\begin{corollary}\label{cor:chow-positive-formula}
The Chow polynomial of a connected lattice path matroid $M$ of rank $r$ is
    \begin{equation*}
        \Chow_M(x)
        = \sum_{S} \sum_{\lambda} x^{\des(\lambda)} (1+x)^{r-1-2\des(\lambda)} \,,
    \end{equation*}
    where $\des(\lambda) = |\Des(\lambda)|$ and where the sum ranges over all connected snakes~$S$ contained in~$M$ and all $(M,S)$-surviving sequences $\lambda \in \Lab(S)$ for which $i \in\{0\} \cup \Des(\lambda)$ implies $i+1 \notin \Des(\lambda)$.
\end{corollary}
\begin{proof}
    By~\cite[Theorem 1.1]{stump2024chowaugmentedchowpolynomials},
    the Chow polynomial is determined by the flag $h$-vector as
    \[
        \Chow_M(x) = \sum_{T} \beta_{\cL}(T) \cdot x^{|T|} \, (1+x)^{r-1-2|T|} \,, 
    \]
    where the sum ranges over subsets $T\subseteq\{2,\dots, r-1\}$ with no two consecutive elements.
    The proof follows with \Cref{cor:flag-h-LPM-via-snakes}.
\end{proof}

\subsection{Combinatorial description of \texorpdfstring{$(M,S)$}{(M,S)}-surviving sequences}\label{subsec:intro-MS-surviving}

Let $S = M[U_S,L_S]$ be a connected snake of rank~$r$ on $[n]$ with upper path $U_S=u_1 u_2\cdots u_n$.
Since $S$ is connected, we have $u_1=\lp{N}$ and $u_n=\lp{E}$, and the lower path is obtained by exchanging these two steps, $L_S = \lp{E}\, u_2\cdots u_{n-1}\, \lp{N}$.
Following \cite{Knauer2018TutteLatticePath}, we recall the construction of a graph $G_S$ having~$S$ as its graphic matroid.
It has $r+1$ vertices $v_0,v_1,\dots,v_r$, where $v_0$ is called the \Dfn{root}, and $n$~edges labeled by $1,\dots,n$, one for each step of~$U_S$.
Let $1=y_1<y_2<\dots<y_r$ be the positions of the north steps in~$U_S$.
We first construct the \Dfn{initial path} of~$G_S$
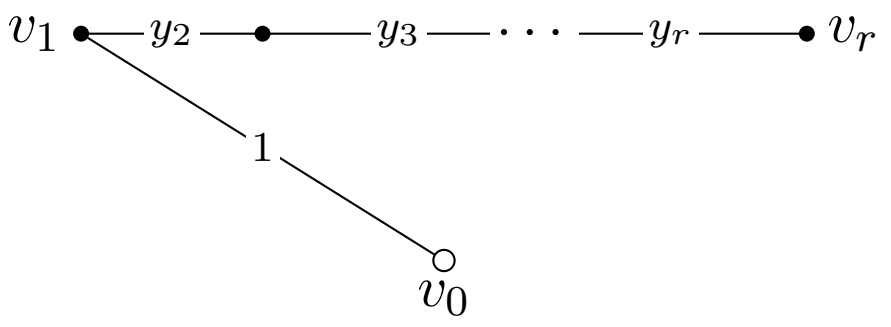
\[
\begin{tikzpicture}[lbl/.style={fill=white, inner sep=1pt, font=\scriptsize}]
    \coordinate (v0) at (0,0);
    \coordinate (v1) at (-2.4,1.5);
    \coordinate (v2) at (-1.2,1.5);
    \coordinate (vd) at (0.6,1.5);
    \coordinate (vr) at (2.4,1.5);
    \draw (v0) -- node[lbl] {$1$} (v1);
    \draw (v1) -- node[lbl] {$y_2$} (v2);
    \draw (v2) -- node[lbl] {$y_3$} (vd);
    \draw (vd) -- node[lbl] {$y_r$} (vr);
    \node[fill=white, inner sep=1pt] at (vd) {$\cdots$};
    \fill (v1) circle (1.5pt);
    \fill (v2) circle (1.5pt);
    \fill (vr) circle (1.5pt);
    \filldraw[fill=white] (v0) circle (2pt);
    \node[left] at (v1) {$v_1$};
    \node[right] at (vr) {$v_r$};
    \node[below] at (v0) {$v_0$};
\end{tikzpicture}
\]
in which the edge between~$v_{k-1}$ and~$v_k$ is labeled by~$y_k$ for $k=1,\dots,r$, see the middle of \Cref{fig:graph-construction}.
In other words, the $k$-th north step of~$U_S$, going from height~$k-1$ to height~$k$, becomes the edge between $v_{k-1}$ and $v_k$.
The remaining $n-r$ edges connect vertices among $v_1,\dots,v_r$ to the root:
an east step~$u_i$ at height~$k$, that is, an index~$i$ with $\#\{y_j\mid y_j<i\}=k$, becomes an edge between~$v_k$ and the root~$v_0$.
Since $u_n = \lp{E}$, we in particular have added the edge between~$v_r$ and the root~$v_0$ labeled by~$n$, thus the graph $G_S$ is the cycle $\{y_1<\dots<y_{r} < n\}$ with further edges incident to the root.

\begin{example}\label{ex:graph-construction}
    Let $S=\LPM{NENENE}{EENENN}$ be the connected snake of rank~$3$ on~$[6]$ shown in \Cref{fig:graph-construction}.
    The north steps of~$U_S$ are at positions $(y_1,y_2,y_3)=(1,3,5)$, so the initial path has the edges $1,3,5$, where the vertex~$v_k$ corresponds to height~$k$.
    The east steps $2$, $4$ and $6$ of~$U_S$ at heights~$1$, $2$ and $3$ become the edges from $v_1$, $v_2$ and $v_3$ to the root.
\end{example}
\begin{figure}[H]
    \centering
    \begin{tikzpicture}[scale=1, baseline=(current bounding box.center)]
        \drawcells[snakestyle]{0/0, 1/0, 1/1, 2/1, 2/2}
        \drawpath{E,E,N,E,N,N}
        \drawpath*[very thick]{N,E,N,E,N,E}
        \foreach \k in {0,...,3}
            \node[right, font=\scriptsize] at (3.15,\k) {$v_\k$};
    \end{tikzpicture}
    \qquad
    \begin{tikzpicture}[baseline=(current bounding box.center)]
        \drawgraph{3}{1/0/1/0/e, 3/1/2/0/e, 5/2/3/0/e}
        \node[left] at (v1) {$v_1$};
        \node[above] at (v2) {$v_2$};
        \node[right] at (v3) {$v_3$};
        \node[below] at (v0) {$v_0$};
    \end{tikzpicture}
    \qquad
    \begin{tikzpicture}[baseline=(current bounding box.center)]
        \drawgraph{3}{1/0/1/0/e, 3/1/2/0/e, 5/2/3/0/e, 6/3/0/0/e, 2/1/0/30/e, 4/2/0/0/e}
        \node[left] at (v1) {$v_1$};
        \node[above] at (v2) {$v_2$};
        \node[right] at (v3) {$v_3$};
        \node[below] at (v0) {$v_0$};
    \end{tikzpicture}
    \caption{Construction of $G_S$ in \Cref{ex:graph-construction}.}
    \label{fig:graph-construction}
\end{figure}

\medskip
Now let the connected snake~$S$ be contained in some connected LPM~$M=M[U,L]$.
A \Dfn{valley} of~$U_S$ is a position~$k$ with $u_{k-1}u_k=\lp{EN}$, and we call it a \Dfn{low valley} if the path obtained from~$U_S$ by changing this $u_{k-1}u_k=\lp{EN}$ into $\lp{NE}$ is weakly below the upper path~$U$ of~$M$.
Define $G_S^M$ as the graph $G_S$ in which all low valleys are marked (by circles in \textcolor{red}{red}).
Observe that a low valley $k=y_j$ is an edge of the initial path not incident to the root, while the preceding edge~$k-1$ is incident to the root and shares the vertex~$v_{j-1}$ with~$k$.

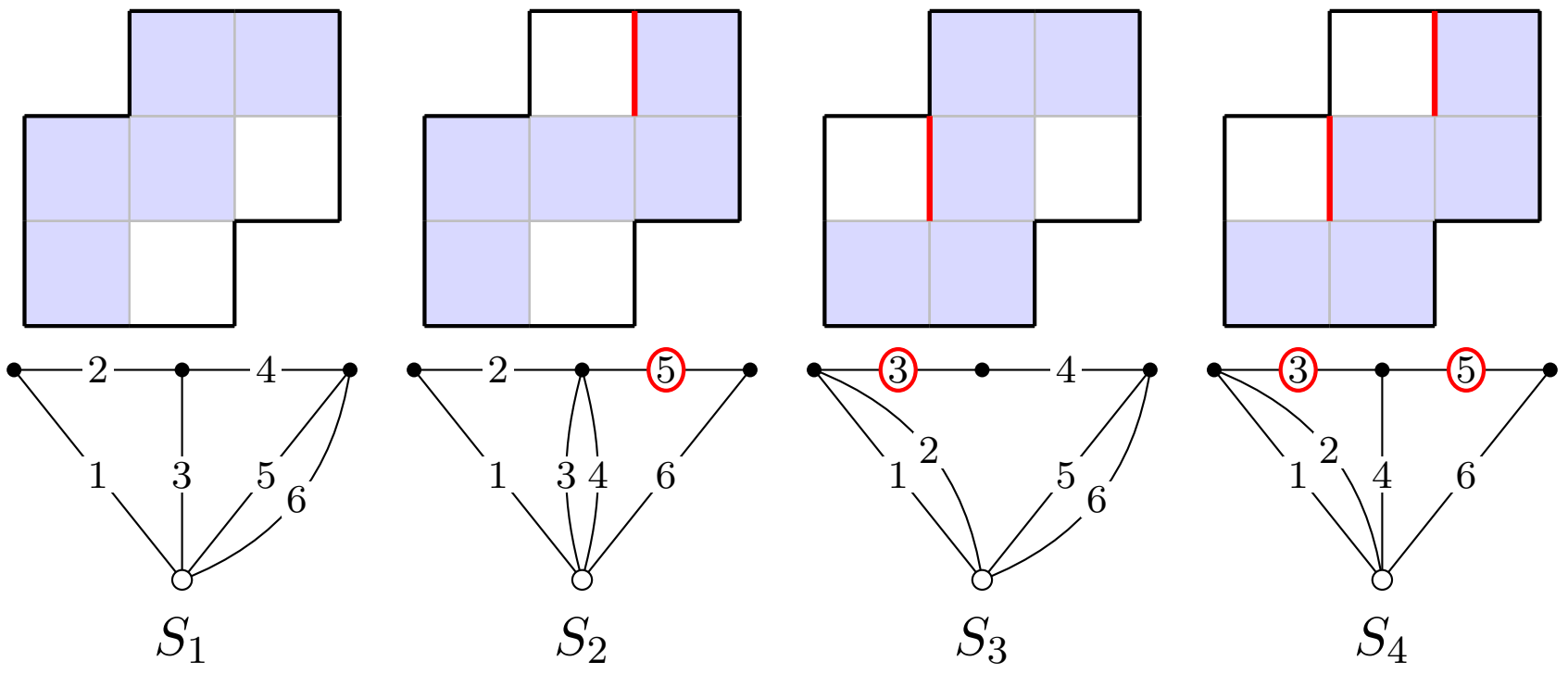
\begin{figure}[H]
    \centering\begin{tabular}{cccc}
         \begin{tikzpicture}[scale=0.75]
             \drawcells[snakestyle]{0/0, 0/1, 1/1, 1/2, 2/2}
             \drawM
         \end{tikzpicture} &
         \begin{tikzpicture}[scale=0.75]
             \drawcells[snakestyle]{0/0, 0/1, 1/1, 2/1, 2/2}
             \drawM
             \draw[red, very thick] (2,2) -- (2,3);
         \end{tikzpicture} &
         \begin{tikzpicture}[scale=0.75]
             \drawcells[snakestyle]{0/0, 1/0, 1/1, 1/2, 2/2}
             \drawM
             \draw[red, very thick] (1,1) -- (1,2);
         \end{tikzpicture} &
         \begin{tikzpicture}[scale=0.75]
             \drawcells[snakestyle]{0/0, 1/0, 1/1, 2/1, 2/2}
             \drawM
             \draw[red, very thick] (1,1) -- (1,2) (2,2) -- (2,3);
         \end{tikzpicture} \\
         \begin{tikzpicture}
             \drawgraph{3}{1/0/1/0/e, 2/1/2/0/e, 4/2/3/0/e, 5/3/0/0/e, 3/2/0/0/e, 6/3/0/30/e}
         \end{tikzpicture} &
         \begin{tikzpicture}
             \drawgraph{3}{1/0/1/0/e, 2/1/2/0/e, 5/2/3/0/lv, 6/3/0/0/e, 3/2/0/-15/e, 4/2/0/15/e}
         \end{tikzpicture} &
         \begin{tikzpicture}
             \drawgraph{3}{1/0/1/0/e, 3/1/2/0/lv, 4/2/3/0/e, 5/3/0/0/e, 2/1/0/30/e, 6/3/0/30/e}
         \end{tikzpicture} &
         \begin{tikzpicture}
             \drawgraph{3}{1/0/1/0/e, 3/1/2/0/lv, 5/2/3/0/lv, 6/3/0/0/e, 2/1/0/30/e, 4/2/0/0/e}
         \end{tikzpicture} \\
         $S_1$ & $S_2$ & $S_3$ & $S_4$
    \end{tabular}
    \caption{The four connected snakes $S_1,\dots,S_4$ contained in $M$ (top) and the graphs $G^M_{S_i}$ (bottom). Low valleys are marked in red, the root is white.}
    \label{fig:running-snakes}
\end{figure}

By \cite[Theorem 3.2]{Knauer2018TutteLatticePath}, the graphic matroid of $G_S$ is $S$.
We therefore identify the ground set $[n]$ of $S$ with the edge set of $G_S^M$, and subgraphs of $G_S^M$ with their edge sets.
Under this identification, the bases of $S$ are the spanning trees of $G_S^M$, the circuits of $S$ are its cycles, and the fundamental circuit of an edge $x\notin T$ with respect to a spanning tree~$T$ is the \Dfn{fundamental cycle} $C(x,T)$, the unique cycle in $T\cup\{x\}$.
In particular, for an $S$-sequence $\lambda=(\lambda_1,\dots,\lambda_r)$, the set $T_\lambda=\{\lambda_1,\dots,\lambda_r\}$ is a spanning tree of $G_S^M$, and \Cref{prop:LabM-general} reads as follows:
an integer sequence~$\lambda$ is an $S$-sequence if and only if $T_\lambda$ is a spanning tree of $G_S^M$ and, for every edge $x\notin T_\lambda$, the subsequence $\lambda|_{C(x,T_\lambda)}$ ends with an integer smaller than~$x$.

\medskip
The $(M,S)$-surviving sequences are characterized by a condition of the same type, imposed at the low valleys.
Recall that the edges of $G_S^M$ not incident to the root are exactly the edges $y_2,\dots,y_r$ of the initial path, and that every low valley is one of them.
For a set $C\subseteq[n]$ and $y\in[n]$, we write $C_{\geq y}=\{z\in C\mid z\geq y\}$.

\begin{theorem}\label{thm:intro-LPM-surviving-seq}
    Let $S$ be a connected snake contained in $M$.
    An $S$-sequence~$\lambda$ is $(M,S)$-surviving if and only if, for every low valley $y\in T_\lambda$, there is an edge $x\notin T_\lambda$ not incident to the root with $y\in C(x,T_\lambda)$ such that $\lambda|_{C(x,T_\lambda)_{\geq y}}$ ends with an integer larger than~$x$.
\end{theorem}

The edge~$x$ in \Cref{thm:intro-LPM-surviving-seq} is unique.
We show in \Cref{lemma:lowvalley-in-fundamental-cycle} that it is the smallest edge $x\geq y$ that is neither in~$T_\lambda$ nor incident to the root.

\begin{example}\label{ex:running-intro}
    Consider the snake $S_4=\LPM{NENENE}{EENENN}$ in \Cref{fig:running-snakes} with low valleys $3$ and $5$, and the spanning tree $T=\{1,3,6\}$ of $G_{S_4}^M$ shown on the left of \Cref{fig:running-intro-trees}.
    The $S_4$-sequences $\lambda$ with $T_\lambda=T$ are $(1,6,3)$, $(3,6,1)$, $(6,1,3)$ and $(6,3,1)$.
    The low valley $5 \notin T$ does not impose any condition.
    For $y=3$, we have $x=5$ with $C(5,T)=\{1,3,5,6\}$, so we look at the entries $3$ and $6$.
    In $(3,6,1)$ the last one is $6>5$, hence this sequence is $(M,S_4)$-surviving.
    In the other three sequences the last one is $3<5$, so these are not.
    The spanning tree $\{1,4,6\}$ on the right of \Cref{fig:running-intro-trees} carries the unique $S_4$-sequence $(6,4,1)$, which is $(M,S_4)$-surviving as it contains no low valley.
    One checks similarly that the remaining spanning trees of $G_{S_4}^M$ carry no $(M,S_4)$-surviving sequences, so that $(3,6,1)$ and $(6,4,1)$ are the only $(M,S_4)$-surviving sequences.
\end{example}
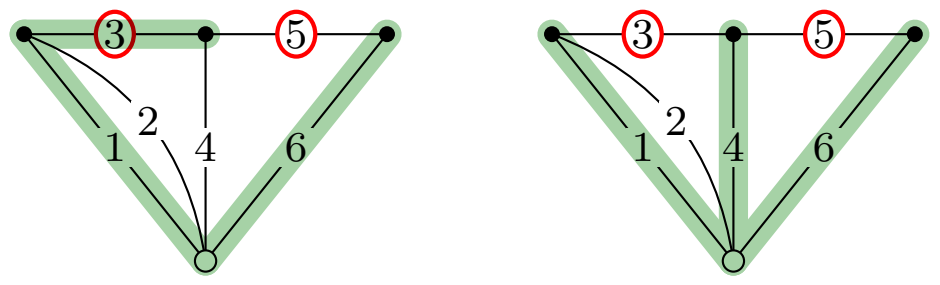
\begin{figure}[H]
    \centering
    \begin{tikzpicture}
        \drawgraph{3}{1/0/1/0/T, 3/1/2/0/lvT, 5/2/3/0/lv, 6/3/0/0/T, 2/1/0/30/e, 4/2/0/0/e}
    \end{tikzpicture}
    \qquad
    \begin{tikzpicture}
        \drawgraph{3}{1/0/1/0/T, 3/1/2/0/lv, 5/2/3/0/lv, 6/3/0/0/T, 2/1/0/30/e, 4/2/0/0/T}
    \end{tikzpicture}
    \caption{The graph $G_{S_4}^M$ of the running example with the spanning trees $T=\{1,3,6\}$ (left) and $T=\{1,4,6\}$ (right) shaded thick in green.
    Low valleys are red.}
    \label{fig:running-intro-trees}
\end{figure}

\subsection{An application to Schubert matroids}\label{subsec:intro-schubert}
A Schubert matroid~$M$ on $[n]$ of rank~$r$ is a lattice path matroid with trivial upper path $\lp{N}^r\lp{E}^{n-r}$.
For a connected snake~$S$ in~$M$, every valley of the upper path of $S$ is a low valley in $(M,S)$. 
We thus call $(M,S)$-surviving sequences \Dfn{$S$-surviving}.
Recall that a connected snake is uniquely defined by its upper path.
For a composition $c=(c_1,\dots , c_\ell) \vDash r-1 \geq 1$,
define two connected snakes of rank~$r$,
$S_c$ with upper path $\lp{N}\,\lp{N}^{c_1}\, \lp{E}\,\lp{N}^{c_2}\, \lp{E}\dots \lp{N}^{c_\ell}\,\lp{E}$ and $S_{c}'$ with upper path $\lp{NE}\,\lp{N}^{c_1}\, \lp{E}\,\lp{N}^{c_2}\, \lp{E}\dots \lp{N}^{c_\ell}\,\lp{E}$\,.

\begin{figure}[H]
    \centering
    \tikzset{
        composition brace/.style={
            decorate, decoration={brace, amplitude=3pt}, thin
        },
        composition label/.style={font=\scriptsize, inner sep=2pt}
    }
    \begin{tikzpicture}[scale=0.65, baseline=(current bounding box.center)]
        \drawcells[snakestyle]{0/0,0/1,0/2,1/2,1/3,2/3,2/4,2/5}
        \drawpath{E,N,N,E,N,E,N,N,N}
        \drawpath{N,N,N,E,N,E,N,N,E}
        \draw[red, very thick]
            (1,3) -- (1,4)
            (2,4) -- (2,5);

        \node[composition label, left] at (-0.15,0.5) {$1$};
        \draw[composition brace] (-0.18,1) -- (-0.18,3)
            node[midway, left=3pt, composition label] {$c_1=2$};
        \draw[composition brace] (0.82,3) -- (0.82,4)
            node[midway, left=3pt, composition label] {$c_2=1$};
        \draw[composition brace] (1.82,4) -- (1.82,6)
            node[midway, left=3pt, composition label] {$c_3=2$};

    \end{tikzpicture}
    \qquad
    \begin{tikzpicture}[scale=0.65, baseline=(current bounding box.center)]
        \drawcells[snakestyle]{0/0,1/0,1/1,1/2,2/2,2/3,3/3,3/4,3/5}
        \drawpath{E,E,N,N,E,N,E,N,N,N}
        \drawpath{N,E,N,N,E,N,E,N,N,E}
        \draw[red, very thick]
            (1,1) -- (1,2)
            (2,3) -- (2,4)
            (3,4) -- (3,5);

        \node[composition label, left] at (-0.15,0.5) {$1$};
        \draw[composition brace] (0.82,1) -- (0.82,3)
            node[midway, left=3pt, composition label] {$c_1=2$};
        \draw[composition brace] (1.82,3) -- (1.82,4)
            node[midway, left=3pt, composition label] {$c_2=1$};
        \draw[composition brace] (2.82,4) -- (2.82,6)
            node[midway, left=3pt, composition label] {$c_3=2$};

    \end{tikzpicture}
    \caption{The snakes $S_c$ and $S'_c$ for the composition $c = (2,1,2) \vDash 5$.}
    \label{fig:snakes-Sc}
\end{figure}
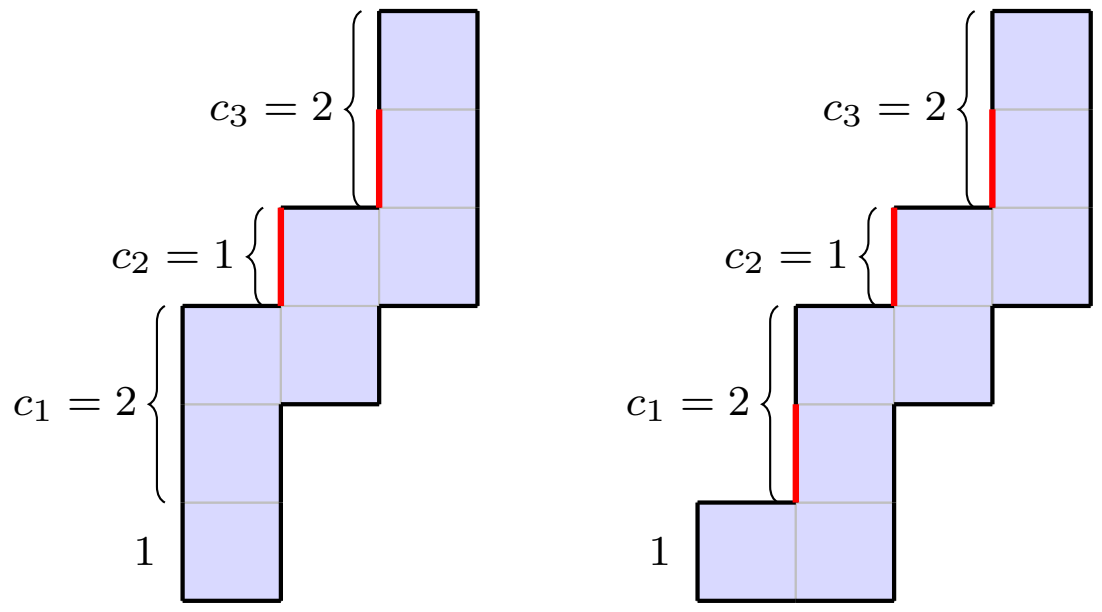

\begin{theorem}\label{thm:ab-schubert}
    The $\aaa\bbb$-index of a connected Schubert matroid $M$ of rank $r\geq 2$ with lower path $L$ equals
    \begin{equation*}
        \Psi_M(\aaa,\bbb)
        = \sum_{c} \Big(f_{c}(L) \sum_{\lambda\  S_{c}\text{-surviving}} \operatorname{m}_{\Des(\lambda)}(\aaa,\bbb)
        + f'_{c}(L) \sum_{\lambda\ S_{c}'\text{-surviving}} \operatorname{m}_{\Des(\lambda)}(\aaa,\bbb)\Big)\,,
    \end{equation*}
    where the sum ranges over all compositions $c \vDash r-1$ and where $f_{c}(L)$ and $f_{c}'(L)$ are integers explicitly defined in \Cref{sec:applications}.
\end{theorem}

\subsection{Outline}
In~\Cref{sec:preliminaries}, we recall background on matroids, matroid subdivisions, and valuative invariants.
In \Cref{sec:M-sequences}, we characterize $M$-sequences and provide a proof of~\Cref{prop:LabM-general}.
In \Cref{sec:snake-subdivision}, we then focus on lattice path matroids and the snake subdivision.
We regroup terms of the $\aaa\bbb$-index into polynomials depending on connected snakes contained in $M$.
Using an interpretation of $S$-sequences for connected snakes $S$ in terms of the graph $G_S$, we then show in \Cref{subsec:MS-surviving-seq} that each of these polynomials admits a positive combinatorial formula in terms of $(M,S)$-surviving sequences, and prove \Cref{thm:intro-ab-via-MS} in \Cref{subsec:MS-surviving-seq}.
Finally, in \Cref{sec:applications}, we apply our results from \Cref{sec:snake-subdivision} to the $\aaa\bbb$-index of Schubert matroids and uniform matroids.

\subsection*{Usage of AI}
\Cref{prop:LabM-general} was suggested by Claude Opus 5.5 as a generalization of the case of graphic matroids. The proof was provided independently by the authors. Claude Opus 5.5 and ChatGPT checked all examples, the proofs and pointed out gaps, among them in the proof of \Cref{thm:intro-LPM-surviving-seq} and a circular argument in \Cref{sec:M-sequences}.

The same tools were used to improve the language, to check inconsistency of notation and cross-references, and to generate the figures from specifications provided by the authors. In this process, they suggested alternative formulations of \Cref{lemma:lowvalley-in-fundamental-cycle} and \Cref{cor:Lab-Sb-condition}, leaving the mathematics unchanged.

Claude Opus 5.5 partially generated the SageMath code to compute $M$-sequences for small lattice path matroids. The results in \Cref{subsec:MS-surviving-seq} were conjectured using SageMath code written by Claude Opus 5.5.

\subsection*{Acknowledgment}
This project started during the \emph{Dive into Research} in Bochum 2025, funded by the DFG priority program SPP 2458 Combinatorial Synergies.
The authors thank Timo Hauser for the collaboration during the program.
The authors also thank Luis Ferroni for an insightful discussion that encouraged them to revisit the project.
CG was supported by ``the Combinatorial Polytope Theory'' (COMPOTE) project (DFG 541393733).
EH was partially supported by the Italian Ministry of University and Research (MUR) under the grant FIS 2, project ``Chow Ring of Toric Arrangements'' (FIS-2023-00808, CUP J53C25001900001).

\section{Preliminaries}\label{sec:preliminaries}
We refer the reader to \cite{Oxley2011} for background on matroid theory and to \cite{BoninDeMier2006} for background on lattice path matroids.
Here, we only recall the definition and properties of the closure operator, of circuits, and of the minimal labeling, which we will mainly use in \Cref{sec:M-sequences}.
We also provide background on matroid subdivisions and valuative invariants, following an overview from \cite[Section 2.2]{ferroni2024valuative}.

\subsection{Closure and circuits}\label{subsec:prelim-matroid-basics}
Let $M$ be a matroid on the ground set $E$, with the usual notions of bases, independence, and rank function $\rk:2^E \to \NN_0$.
The \Dfn{closure} of a subset $A$ is defined as $\closure(A) = \{ x\in E \mid \rk(A\cup\{x\}) = \rk(A) \}$.
A subset $F\subseteq E$ is a \Dfn{flat} if $\closure(F) = F$.
The \Dfn{lattice of flats}~$\cL(M)$ is obtained by ordering the flats by inclusion.
A \Dfn{circuit} is a dependent subset~$C$ for which $C\setminus\{x\}$ is independent for every $x\in C$.
We recall a result from Oxley on sets that contain a (unique) circuit.

\begin{lemma}[{\cite[Corollary 1.2.6]{Oxley2011}}]\label{lemma:Oxley-fund-circ}
    For any basis~$B$ and any $x\in E\setminus B$, there exists a unique circuit $C(x,B) \subseteq B\cup \{x\}$.
\end{lemma}

We obviously have $x \in C(x,B)$, and call this circuit the \Dfn{fundamental circuit of $x$ (with respect to $B$)}.
A \Dfn{loop}~$x\in E$ is an element not contained in any basis of $M$.
Therefore, $x$ is a loop if and only if $C(x,B) = \{x\}$ (with respect to any basis).
The closure of a set can be described in terms of circuits.
\begin{lemma}[{\cite[Proposition 1.4.11]{Oxley2011}}]\label{lemma:Oxley-closure-circ}
    For $A\subseteq E$, we have
    \[
        \closure(A) = A\cup\{x\in E \setminus A \mid \text{there exists a circuit } C \subseteq A \cup \{x\} \text{ with } x\in C\}\,.
    \]
\end{lemma}

\subsection{Flag \texorpdfstring{$f$}{f}-vectors and flag \texorpdfstring{$h$}{h}-vectors}\label{subsec:flag-vectors}
We use the terminology for posets from~\cite{StanleyEC1}.

Let $P$ be a finite graded poset of rank $r$ with unique minimal element $\zero$ and unique maximal element $\one$.
The \Dfn{flag $f$-vector} of $P$ is the map $\alpha_P: 2^{\{0,\dots , r\}}\to \ZZ$ that sends a subset $S\subseteq \{0,\dots , r\}$ to the number of maximal chains in the rank-selected subposet $P_S = \{x\in P \mid \rk(x) \in S\}\cup\{\zero, \one\}$.
The \Dfn{flag $h$-vector} $\beta_P: 2^{\{0,\dots , r\}}\to \ZZ$ is then defined as 
\begin{equation*}
    \beta_P(T) = \sum_{S\subseteq T} (-1)^{|T\setminus S|} \, \alpha_P(S)\,,
    \quad \text{for } T\subseteq \{0,\dots , r\}\,.
\end{equation*}

An edge-labeling of $P$ is a map $\lambda$ from the set of cover relations (or edges) of $P$ to $\ZZ$. 
For a saturated chain $\mathfrak{m}=\{m_0 \prec m_1 \prec \cdots \prec m_k\}$, we write $\lambda(\mathfrak{m})=(\lambda_1,\dots,\lambda_k)$ for the sequence of labels $\lambda_i = \lambda(m_{i-1}\prec m_i)$ along $\mathfrak{m}$.
We call $\lambda$ an \Dfn{$R$-labeling},
if for any two elements~$x<y$ in~$P$ there exists a unique saturated chain~$\mathfrak{m}=\{ x= m_0 \prec m_1 \prec \cdots \prec m_k = y \}$ such that~$\lambda(\mathfrak{m})$ has increasing edge-labels,~$\lambda_{i} < \lambda_{i+1}$ for all $i$.
In this case, the flag $h$-vector is determined by the descents of the edge-labels along maximal chains in $P$:
\[
    \beta_P(S) = \# \{ \mathfrak{m} \text{ maximal chain in }P \mid \Des(\lambda(\mathfrak{m})) = S \}\,,
    \quad
    \text{for } S\subseteq\{0,\dots,r\}\,.
\]
If $P=\cL(M)$ is the lattice of flats of a matroid $M$, fix an order on the ground set.
Björner~\cite[§7.6]{Bjorner1992} shows that the \Dfn{minimal labeling}\footnote{We are stating a mild variation of the minimal labeling.} $\lambda$, which sends a cover $F\prec G$ in $\cL(M)$ to
\begin{equation*}
    \lambda(F\prec G) = \min\{i \mid F \vee \closure(i) = G \}\,,
\end{equation*}
is an $R$-labeling of $\cL(M)$.

\subsection{Matroid subdivision}
The (\emph{base}) \emph{polytope} of the matroid~$M$ is
\[
    \cP(M) 
    = \operatorname{conv}\biggl\{\sum_{i\in B} e_i \ \biggm| \ B\in \mathcal{B}(M)\biggr\} \subseteq \mathbb{R}^E \,,
\]
where $e_i$ is the $i$-th canonical vector in $\RR^E$ and $\mathcal{B}(M)$ is the set of bases.
The dimension of $\cP(M)$ is $n-c(M)$ where~$c(M)$ is the number of connected components of $M$ \cite[Proposition 2.4]{FeichtnerSturmfels2005}.
Following \cite[Definition 2.11]{ferroni2024valuative},
a \Dfn{matroid subdivision} of $\mathcal P(M)$ is a polyhedral subdivision $\mathcal{S}$ into base polytopes of matroids $M_1,\dots , M_s$ on the same ground set~$E$:
\begin{enumerate}
    \item $\cP(M)=\cup_{i=1}^s \cP(M_i)$;
    \item for each $i$, any face of $\cP(M_i)$ belongs to $\mathcal{S}$; and
    \item $\cP(M_i)\cap \cP(M_j)$ is a common face of $\cP(M_i)$ and $\cP(M_j)$.
\end{enumerate}
We denote the set of \Dfn{internal} faces by $\mathcal{S}^{\text{int}} = \{\cP(M_i)\in \mathcal{S} \mid \cP(M_i) \not\subseteq \partial\cP \}$, where $\partial\cP$ is the boundary of $\cP(M)$.

\subsection{Valuative invariants}
We use two notations for collections of matroid polytopes:
\begin{equation*}
    \mathfrak{M}_E  =  \{ \cP(M) \mid M \text{ matroid on } E\}\,, \qquad
    \mathfrak{M} =  \bigsqcup_{E\subset\NN \text{ finite}}\mathfrak{M}_E \,.
\end{equation*}
A map $f$ that sends matroid polytopes into an abelian group $A$, $f:\mathfrak{M} \to A$, is called a \Dfn{valuation} if for any matroid subdivision $\mathcal{S}$ of $\cP\in \mathfrak{M}$, we have
\begin{equation*}
    f(\cP) = \sum_{\mathcal{Q}\in \cS^{\text{int}}} (-1)^{\dim\cP -\dim \mathcal{Q} } \, f(\mathcal{Q}) \,.
\end{equation*}
The map $f$ is called an \Dfn{invariant}, if $f(\cP(M)) = f(\cP(N))$ whenever $M\cong N$ are isomorphic matroids.
If $f$ is both a valuation and an invariant,
we call $f$ a \Dfn{valuative invariant}. 
In this case, we understand $f$ as a map on matroids with $f(M) = f(\cP(M))$.

The map that associates to each matroid $M$ the flag $f$-vector $\alpha_\cL$ is a valuative invariant by \cite[Theorem 8.2]{ferroni2024valuative}.
Since each flag $h$-vector entry is a fixed linear combination of flag $f$-vector entries, we immediately obtain the following corollary.

\begin{corollary}\label{cor:ab-index-valuative}
    The map $M \mapsto \Psi_M(\aaa,\bbb)\in \ZZab$ is a valuative invariant.
\end{corollary}

\section{\texorpdfstring{$M$}{M}-sequences}\label{sec:M-sequences}
In this section, we prove \Cref{prop:LabM-general}.
An $M$-sequence $\lambda$, $\lambda \in \Lab(M)$, determines a chain~$\cF=\{F_0 \prec F_1 \prec \dots \prec F_r\}$ with
\[
    F_i = \closure(\{\lambda_1,\dots , \lambda_i\}) \quad \text{for } 0\leq i \leq r\,,
\]
for which $\lambda(\cF) = \lambda$, since $\lambda_i\in F_i\setminus F_{i-1}$ and $\rk(F_i)=i$ for all $i$.
In particular, the set $B_\lambda = \{\lambda_1,\dots , \lambda_r\}$ is a basis of $M$.
Indeed, it is well known that $M$-sequences form (\emph{NBC}-)bases of $M$, see \Cref{rem:NBC}.
The first condition in~\Cref{prop:LabM-general} is thus a necessary condition.
The second condition comes from the following lemma.

\begin{lemma}\label{lma:circuits-and-M-seq}
    Let $\lambda=(\lambda_1,\dots , \lambda_r)$ be a sequence such that $B_\lambda = \{\lambda_1,\dots , \lambda_r\}$ is a basis of $M$, and set
    $F_i = \closure(\{\lambda_1,\dots , \lambda_i\})$ for $0\leq i\leq r$.
    Let $x\notin B_\lambda$ be a non-loop of $M$ with fundamental circuit $C_x=C(x,B_\lambda)$.
    Then, $\lambda|_{C_x}$ is non-empty, and $x\in F_i\setminus F_{i-1}$, where $\lambda_i$ is the last entry in $\lambda|_{C_x}$.
\end{lemma}
\begin{proof}
    Since~$x$ is not a loop, we have $C_x\neq\{x\}$,
    and hence~$\lambda|_{C_x}$ is non-empty, ending with~$\lambda_i$ for some $i\geq 1$.
    Let $0\leq j\leq r$.
    Every circuit $C$ with $x\in C\subseteq\{\lambda_1,\dots,\lambda_j,x\}$ is contained in~$B_\lambda\cup\{x\}$,
    and hence $C=C_x$ by \Cref{lemma:Oxley-fund-circ}.
    Since $x\notin B_\lambda$,
    \Cref{lemma:Oxley-closure-circ} gives that~$x\in F_j$ if and only if~$C_x\setminus\{x\}\subseteq\{\lambda_1,\dots,\lambda_j\}$.
    This holds if and only if $j\geq i$, and therefore~$x\in F_i\setminus F_{i-1}$.
\end{proof}

\begin{proof}[Proof of \Cref{prop:LabM-general}]
    We have already seen that an~$M$-sequence $\lambda=(\lambda_1,\dots , \lambda_r)$ satisfies the first condition, $B_\lambda$ is a basis.
    Moreover, $\lambda$, with $B_\lambda$ a basis, determines the maximal chain $\cF$ for which $\lambda(\cF)=\lambda$ by setting $F_i = \closure(\{\lambda_1,\dots , \lambda_i\})$ for $0\leq i\leq r$.
    By definition, $\lambda$ is an $M$-sequence if and only if
    \[
        \lambda_i = \min (F_i\setminus F_{i-1}) \quad \text{for all } 1\leq i\leq r\,.
    \]
    It therefore remains to determine the sets $F_i\setminus F_{i-1}$.
    We have $B_\lambda\cap (F_i\setminus F_{i-1}) = \{\lambda_i\}$ and
    every loop of $M$ is contained in $F_0$.
    Every non-loop $x\notin B_\lambda$, by \Cref{lma:circuits-and-M-seq}, lies in $F_i\setminus F_{i-1}$, where $\lambda_i$ is the last entry of $\lambda|_{C(x,B_\lambda)}$.
    Hence, $F_i\setminus F_{i-1}$ consists of $\lambda_i$ and all non-loops $x\notin B_\lambda$ for which $\lambda|_{C(x,B_\lambda)}$ ends with $\lambda_i$.
    Therefore, $\lambda_i=\min(F_i\setminus F_{i-1})$ holds for all $i$ if and only if the second condition in the proposition holds.
\end{proof}

We conclude this section with two remarks on the connection of $M$-sequences to NBC-bases and to valuations.

\begin{remark}[NBC-bases]\label{rem:NBC}
    A basis $B$ of a (loopless) matroid $M$ is an \emph{NBC-basis} if it contains no broken circuit, that is, if there exists no circuit $C$ of $M$ for which~$C\setminus\{\min C\}\subseteq B$,
    see \cite{Rota1964,Bjorner1992}.
    The basis $B$ contains such a broken circuit~$C\setminus\{\min C\}$ if and only if $\min C\notin B$ and $C=C(\min C,B)$ is the fundamental circuit of $\min C$.
    Hence, $B$ is an NBC-basis if and only if
    \[
        x\neq \min C(x,B) \quad \text{for all } x\notin B\,.
    \]
    \Cref{prop:LabM-general} recovers that the sets $B_\lambda$ for $\lambda\in\Lab(M)$ are in bijection with the maximal chains of $\cL$ with decreasing labels, and their number is therefore $\beta_\cL(\{1,\dots , r-1\})$.
\end{remark}

\begin{remark}[$M$-sequences as valuation]
    Denote the set of base polytopes of all matroids on ground set $E$ with total order $\leq$ by $\mathfrak{M}_{E,\leq}$. For an integer sequence $\lambda$ of finite length, define the map $f_\lambda: \mathfrak{M}_{E,\leq} \to \{0,1\}$ by $f_\lambda(\cP(M)) = 1$ if $\lambda\in \Lab(M)$ and $f_\lambda(\cP(M)) = 0$ if $\lambda\notin \Lab(M)$.
    Following the proof of \cite[Theorem 8.2]{ferroni2024valuative}, one can show
    \[
        f_\lambda(\cP(M))
        = \sum_{\cP(Q)\in\mathcal{S}^{\textnormal{int}}}
            (-1)^{\dim\cP(M)-\dim\cP(Q)} \, f_\lambda(\cP(Q)) 
    \]
    for every matroid subdivision $\cS$ of $\cP(M)$ into base polytopes of matroids on the same \emph{ordered} ground set $E$, so that~$\cP(Q)\in \mathfrak{M}_{E,\leq}$ for every internal face $\cP(Q)\in\mathcal{S}^{\textnormal{int}}$.
    Moreover, $f_\lambda$ is invariant under matroid isomorphisms that preserve the order on the ground set.
\end{remark}

\section{The snake subdivision of LPMs}\label{sec:snake-subdivision}
Throughout, let $M=M[U,L]$ be a connected lattice path matroid on $[n]$ of rank~$r$.
Recall from \Cref{subsec:intro-main-results-lpm} that the snake subdivision of $M$ provides \eqref{eq:ab-LPMs},
\begin{equation*}
    \Psi_M(\aaa,\bbb) = \sum_{Q} (-1)^{c(Q)-1} \, \Psi_Q(\aaa,\bbb) \,,
\end{equation*}
where the sum is over all internal faces $Q$ of $M$.
That is, by \cite[Theorem 3.3]{Ferroni2023SchubertMD},
the sum ranges over all matroids~$Q$ such that $Q=M[U_Q,L_Q]$ is a snake for which both~$M[U_Q,L]$ and~$M[U,L_Q]$ are connected.

\subsection{Regrouping terms by low valleys}
We regroup the terms in \eqref{eq:ab-LPMs}
following an idea of Ferroni--Morales--Panova \cite[Section 3.1]{FerroniMoralesPanova2026}.
We present their constructions with a harmless modification\footnote{\cite{FerroniMoralesPanova2026} uses high valleys rather than low peaks.} that is necessary as discussed in \Cref{rem:high-peak-not-boolean}.

\smallskip
A valley of a lattice path $\gamma=\gamma_1\cdots \gamma_n$ is a position $k\geq 2$ such that $\gamma_{k-1}\gamma_k=\lp{EN}$.
For a subset~$A$ of valleys of~$\gamma$,
we write $\gamma^{(A)}$ for the lattice path obtained from $\gamma$ by exchanging $\gamma_{k-1}\gamma_k=\lp{EN}$ by a peak $\lp{NE}$ for all $k\in A$.
For an internal face $Q=M[U_Q,L_Q]$ of some matroid $M=M[U,L]$,
we call a position~$k\geq 2$ a \Dfn{low valley} in $(M,Q)$ if $k$ is a valley of both paths $U_Q$ and $L_Q$, and~$U_Q^{(k)}$ stays weakly below the upper path $U$ of $M$.
The set of low valleys in $(M,Q)$ is denoted by $\lowvalley_M(Q)$.

\smallskip
Let $S=M[U_S,L_S]$ be a connected snake contained in $M$ and let $A\subseteq \lowvalley_M(S)$ be a subset of low valleys in $(M,S)$.
Note that every valley of $U_S$ is a valley of $L_S$,
and we write $S^{(A)}=M[U_S,L_S^{(A)}]$ for the lattice path matroid obtained by exchanging the valleys $A$ in the lower path $L_S$ by peaks.
This gives rise to a bijection from \cite[Lemma 3.4]{FerroniMoralesPanova2026},
\begin{equation}\label{eq:bijection-snakes-lowvalleys}
    \begin{array}{@{}c@{\quad}c@{\quad}c@{}}
        \bigl\{ (S,A) \mid S \text{ connected snake in } M,\,
            A \subseteq \lowvalley_M(S) \bigr\}
        & \longrightarrow
        & \{ \text{internal faces of } M \}
        \\[5pt]
        (S,A)
        & \longmapsto
        & S^{(A)}
    \end{array}
\end{equation}
For a connected snake $S$ contained in $M$, define
\begin{equation}\label{def:ab-incl-exl-for-snake}
    \Psi_{(M,S)}^{\pm}(\aaa,\bbb)
    = \sum_{A\subseteq \lowvalley_M(S)} (-1)^{|A|} \, \Psi_{S^{(A)}}(\aaa,\bbb) \,,
\end{equation}
which regroups the terms for the $\aaa\bbb$-index in \eqref{eq:ab-LPMs} by connected snakes.

\begin{lemma}\label{lma:ab-via-lowvalley}
    The $\aaa\bbb$-index of a connected lattice path matroid~$M$ decomposes as 
    \begin{equation*}
        \Psi_M(\aaa,\bbb)
        = \sum_{S} \Psi_{(M,S)}^{\pm}(\aaa,\bbb) \,,
    \end{equation*}
    where the sum ranges over all connected snakes $S$ contained in $M$.
\end{lemma}
\begin{proof}
    Using the bijection in \eqref{eq:bijection-snakes-lowvalleys}, identify
    every internal face $Q$ of $M$ with the unique connected snake $S$ and the unique subset $A\subseteq\lowvalley_M(S)$ such that $Q=S^{(A)}$.
    We increase the number of components by exactly 1 for every low valley $k\in A$ that is exchanged in the lower path of $S$.
    The proof thus follows with \eqref{eq:ab-LPMs} and $c(S^{(A)})=|A|+1$.
\end{proof}
Note that \Cref{lma:ab-via-lowvalley} is a direct consequence of \cite[Remark 3.8]{FerroniMoralesPanova2026}.

\begin{example}\label{ex:running-ab}
    We continue the running example $M=\LPM{NNENEE}{EENENN}$.
    The nine internal faces are given below with low valleys marked in \textcolor{red}{red} in the lower path.
    
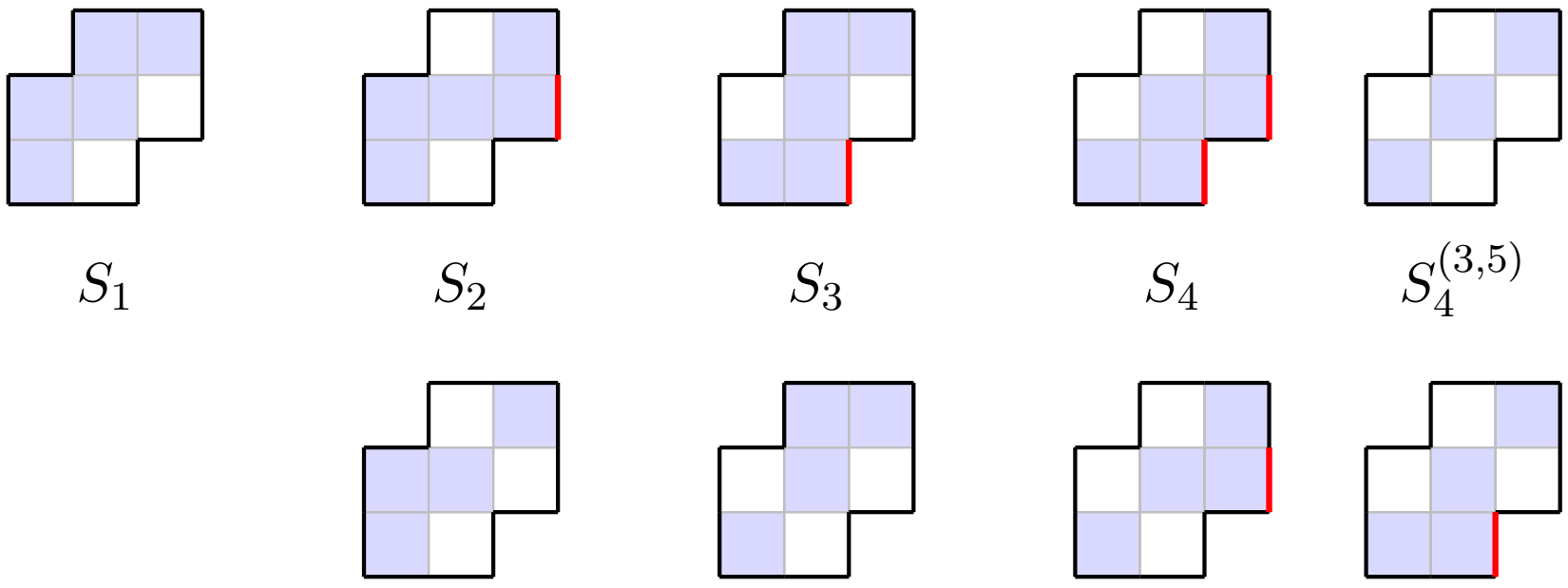
\begin{gather*}
    \begin{tikzpicture}[scale=0.45, baseline=0.55cm]
        \path (-0.2,-2) rectangle (21.2,3.1);
        \begin{scope}[shift={(0,0)}]
            \drawcells[snakestyle]{0/0, 0/1, 1/1, 1/2, 2/2}
            \drawM
            \node[anchor=base] at (1.5,-1.5) {$S_1$};
        \end{scope}
        \begin{scope}[shift={(5.5,0)}]
            \drawcells[snakestyle]{0/0, 0/1, 1/1, 2/1, 2/2}
            \drawM
            \draw[red, very thick] (3,1) -- (3,2);
            \node[anchor=base] at (1.5,-1.5) {$S_2$};
        \end{scope}
        \begin{scope}[shift={(11,0)}]
            \drawcells[snakestyle]{0/0, 1/0, 1/1, 1/2, 2/2}
            \drawM
            \draw[red, very thick] (2,0) -- (2,1);
            \node[anchor=base] at (1.5,-1.5) {$S_3$};
        \end{scope}
        \begin{scope}[shift={(16.5,0)}]
            \drawcells[snakestyle]{0/0, 1/0, 1/1, 2/1, 2/2}
            \drawM
            \draw[red, very thick] (2,0) -- (2,1) (3,1) -- (3,2);
            \node[anchor=base] at (1.5,-1.5) {$S_4$};
        \end{scope}
        \begin{scope}[shift={(21,0)}]
            \drawcells[snakestyle]{0/0, 1/1, 2/2}
            \drawM
            \node[anchor=base] at (1.5,-1.5) {$S_4^{(3,5)}$};
        \end{scope}
    \end{tikzpicture}
    \tag{$+$} \\[4pt]
    \begin{tikzpicture}[scale=0.45, baseline=0.55cm]
        \path (-0.2,-2) rectangle (21.2,3.1);
        \begin{scope}[shift={(5.5,0)}]
            \drawcells[snakestyle]{0/0, 0/1, 1/1, 2/2}
            \drawM
            \node[anchor=base] at (1.5,-1.5) {$S_2^{(5)}$};
        \end{scope}
        \begin{scope}[shift={(11,0)}]
            \drawcells[snakestyle]{0/0, 1/1, 1/2, 2/2}
            \drawM
            \node[anchor=base] at (1.5,-1.5) {$S_3^{(3)}$};
        \end{scope}
        \begin{scope}[shift={(16.5,0)}]
            \drawcells[snakestyle]{0/0, 1/1, 2/1, 2/2}
            \drawM
            \draw[red, very thick] (3,1) -- (3,2);
            \node[anchor=base] at (1.5,-1.5) {$S_4^{(3)}$};
        \end{scope}
        \begin{scope}[shift={(21,0)}]
            \drawcells[snakestyle]{0/0, 1/0, 1/1, 2/2}
            \drawM
            \draw[red, very thick] (2,0) -- (2,1);
            \node[anchor=base] at (1.5,-1.5) {$S_4^{(5)}$};
        \end{scope}
    \end{tikzpicture}
    \tag{$-$}
    \end{gather*}
    The $\aaa\bbb$-indices of the snakes in the first row $(+)$ contribute with a positive sign to~$\Psi_M(\aaa,\bbb)$ and the snakes in the second row $(-)$ with a negative sign.
    To give one example,
    we have 
    \begin{equation*}
    \begin{aligned}[c]
        \Psi_{S_4}(\aaa,\bbb) 
            &= \aaa\aaa + 4\bbb\aaa + 5\aaa\bbb +4\bbb\bbb \\
        \Psi_{S_4^{(3)}}(\aaa,\bbb) 
            &= \aaa\aaa + 3\bbb\aaa + 3\aaa\bbb +2\bbb\bbb \\
    \end{aligned}
    \qquad
    \begin{aligned}[c]
        \Psi_{S_4^{(3,5)}}(\aaa,\bbb) 
            &= \aaa\aaa+2\bbb\aaa+2\aaa\bbb+\bbb\bbb \\
        \Psi_{S_4^{(5)}}(\aaa,\bbb) 
            &= \aaa\aaa + 3\bbb\aaa + 3\aaa\bbb +2\bbb\bbb 
    \end{aligned}
    \end{equation*}
    which gives $\Psi_{(M,S_4)}^{\pm}(\aaa,\bbb)=\aaa\bbb+\bbb\bbb$.
\end{example}

\subsection{\texorpdfstring{$M$}{M}-sequences of snake matroids}\label{subsec:M-seq-snake-matroids}
We show in this section that the coefficients of~$\Psi_{(M,S)}^{\pm}(\aaa,\bbb)\in \ZZab$ are nonnegative by interpreting the formula~\eqref{def:ab-incl-exl-for-snake} as an inclusion-exclusion formula.
This will provide us with a natural definition of $(M,S)$-surviving sequences and a proof of \Cref{thm:intro-ab-via-MS}. 

\medskip
Recall the graphs $G_S$ (and $G_S^M$) from \Cref{subsec:intro-MS-surviving} for a connected snake $S$ contained in~$M$.
The internal faces of $M$ are of the form $S^{(A)}$ for connected snakes $S$ contained in $M$ and subsets $A\subseteq \lowvalley_M(S)$, see \eqref{eq:bijection-snakes-lowvalleys}.
The graph $G_{S,A}$ is obtained from $G_S$ by replacing the endpoint $v_{i-1}$ of $a=y_i$ by the root for every $a\in A$, see \Cref{fig:general-flip}.
Moreover, we mark the low valleys $\lowvalley_M(S^{(A)})$ in $G_{S,A}$ to obtain the marked graph $G_{S,A}^M$.

Together with \cite[Theorem 3.2]{Knauer2018TutteLatticePath}, we obtain the following consequence.
\begin{corollary}\label{cor:graph-after-flips}
    An internal face $S^{(A)}$ is the graphic matroid of the graph $G_{S,A}$.
\end{corollary}
\begin{proof}
    The lower and upper path of the snake $S^{(A)}=M[U_{S^{(A)}},L_{S^{(A)}}]$ uniquely decompose as $U_{S^{(A)}} = U_1\cdots U_m$ and $L_{S^{(A)}} = L_1\cdots L_m$ such that 
    $S_i=M[U_i,L_i]$ is a connected snake on ground set $[n_i]$.
    Let $G_i$ be the graph $G_{S_i}$ with the edge labels shifted by~$n_1+\dots + n_{i-1}$.
    A tree (resp. cycle) in $G_{S,A}$ is then the union of trees (resp. cycles) in these graphs $G_1,\dots, G_m$ of the connected components of $S^{(A)}$.
    Hence, $S^{(A)}$ is the graphic matroid of the graph $G$ obtained from $G_1,\dots, G_m$ by identifying the roots.
    Moreover, $G$ equals $G_{S,A}$ up to a relabeling of the vertices.
    To see this, consider a single low valley $A=\{a\}$.
    The matroid $S^{(a)}$ has two connected components of size $a-1$ and $n+1-a$, respectively, so the claim follows from the construction of $G_{S,\{a\}}$ and induction on the size of $A$.
\end{proof}

\begin{figure}[H]
    \centering
    \tikzset{near top/.style={graph label/.append style={pos=0.25}}}
    \begin{tikzpicture}[x=0.7cm, y=1.3cm]
        \drawgraph{9}{1/1/0/0/near top, {}/1/2/0/dots,
            {}/2/0/0/near top, a/2/3/0/{lv,draw=orange,very thick}, {}/3/4/0/dots,
            {}/4/0/0/near top, b/4/5/0/lv, {}/5/6/0/dots,
            {}/6/0/0/near top, c/6/7/0/{lv,draw=orange,very thick}, {}/7/0/0/e, {}/7/8/0/e,
            {}/8/9/0/dots, n/9/0/0/near top}
    \end{tikzpicture}
    \quad
    \begin{tikzpicture}[x=0.7cm, y=1.3cm]
        \drawgraph{9}{1/1/0/0/near top, {}/1/2/0/dots,
            {}/2/0/0/near top, a/3/0/0/{graph label/.append style={pos=0.45},draw=orange,very thick}, {}/3/4/0/dots,
            {}/4/0/0/near top, b/4/5/0/lv, {}/5/6/0/dots,
            {}/6/0/0/near top, c/7/0/-12/{graph label/.append style={pos=0.45},draw=orange,very thick}, {}/7/0/0/e, {}/7/8/0/e,
            {}/8/9/0/dots, n/9/0/0/near top}
    \end{tikzpicture}
    \caption{The graphs $G^M_{S}$ (left) and $G^M_{S,\{a,c\}}$ (right) for $\{a<b<c\} \subseteq \lowvalley_M(S)$.}
    \label{fig:general-flip}
\end{figure}

As an immediate consequence of the \Cref{cor:graph-after-flips} and \Cref{prop:LabM-general}, we obtain that a sequence $\lambda=(\lambda_1,...,\lambda_r)$ is an $S^{(A)}$-sequence
if and only if
\begin{itemize}
    \item \textbf{Ordered trees:}
        $T_\lambda = \{\lambda_1,\dots ,\lambda_r\}$ is a spanning tree of $G_{S,A}$, and 
    \item \textbf{Fundamental cycles:} for each edge~$x\notin T_\lambda$ with (unique) fundamental cycle $C_x=C(x,T_\lambda)$,
    the sequence~$\lambda|_{C_x}$ ends with an integer smaller than~$x$.
\end{itemize}
Thus, to understand how $S$-sequences differ from $S^{(A)}$-sequences for subsets $A\subseteq \lowvalley_M(S)$, we need to understand how trees and cycles differ.

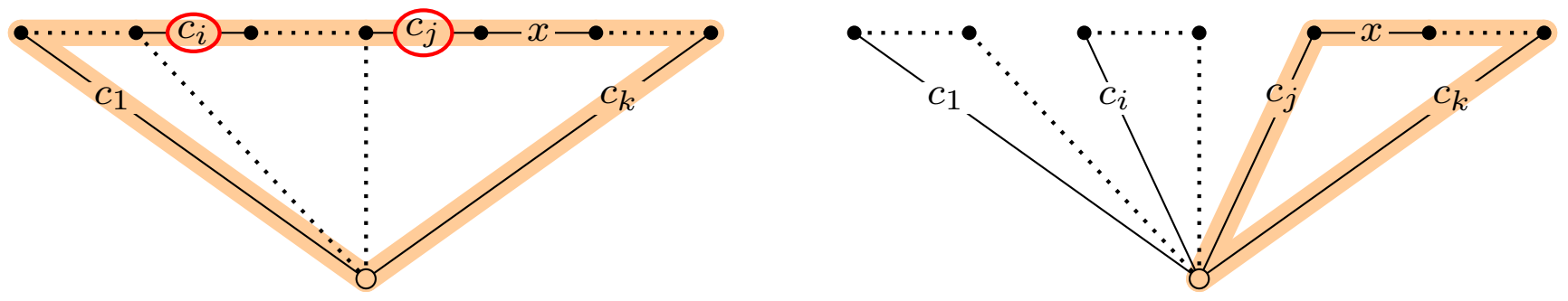
\begin{figure}[H]
    \centering
    \tikzset{near top/.style={graph label/.append style={pos=0.25}}}
    \begin{tikzpicture}[x=0.7cm, y=1.2cm]
        \drawgraph{7}{c_1/1/0/0/{near top,C}, {}/1/2/0/{dots,C}, {}/2/0/0/dots,
            c_i/2/3/0/{lv,C}, {}/3/4/0/{dots,C}, {}/4/0/0/dots,
            c_j/4/5/0/{lv,C}, x/5/6/0/C, {}/6/7/0/{dots,C}, c_k/7/0/0/{near top,C}}
    \end{tikzpicture}
    \qquad
    \begin{tikzpicture}[x=0.7cm, y=1.2cm]
        \drawgraph{7}{c_1/1/0/0/near top, {}/1/2/0/dots, {}/2/0/0/dots,
            c_i/3/0/0/near top, {}/3/4/0/dots, {}/4/0/0/dots,
            c_j/5/0/0/{near top,C}, x/5/6/0/C, {}/6/7/0/{dots,C}, c_k/7/0/0/{near top,C}}
    \end{tikzpicture}
    \caption{Sections of the graphs $G_{S,A}$ (left) and $G_{S,B}$ (right), with $c_i,c_j\in B\setminus A$, and with the fundamental cycles $C_A$ and $C_B$, respectively, of an edge $x\notin T$ shaded in orange.}
    \label{fig:containment-of-cycles}
\end{figure}

\begin{lemma}\label{lemma:containment-of-cycles}
    Let $S$ be a connected snake contained in a lattice path matroid $M$, and let~$A\subseteq B\subseteq \lowvalley_M(S)$.
    \begin{enumerate}[label=(\roman*)]
        \item For $\lambda\in \Lab(S^{(A)})$, the tree $T_\lambda$ is a spanning tree in $G_{S,B}$.
        \item Every tree in $G_{S,B}$ is a tree in $G_{S,A}$.
        \item Let $T$ be a spanning tree in $G_{S,B}$ (and $G_{S,A}$), let $x\notin T$, and denote the fundamental cycle of $x$ in $G_{S,A}$ by $C_A$, and in $G_{S,B}$ by $C_B$.
        Then
        \begin{equation*}
            C_B =
            \begin{cases}
                C_A & \text{if } C_A\cap (B\setminus A)=\emptyset\,, \\
                \{z\in C_A\mid z\geq b\} & \text{otherwise, where } b=\max\big(C_A\cap (B\setminus A)\big)\,.
            \end{cases}
        \end{equation*}
    \end{enumerate}
\end{lemma}
\begin{proof}
    For all parts, we observe that the graphs $G_{S,A}$ and $G_{S,B}$ differ exactly in the edges in $B\setminus A$, which are incident to the root only in $G_{S,B}$.
    Moreover, both the minimal and the maximal edge of any cycle in $G_{S,A}$ or $G_{S,B}$ are incident to the root.

    (i)
    Let $C=\{c_1< \dots < c_k\} \subset T_\lambda$ be a cycle in $G_{S,B}$.
    Since $T_\lambda$ is a tree in $G_{S,A}$ with $c_k$ incident to the root, the edge $c_1\in B\setminus A$ must be the unique low valley in $(M,S^{(A)})$ in $C$.
    Consequently, $c_1$ is not incident to the root in $G_{S,A}$ and $c_1 -1\notin T_\lambda$ is incident to the root.
    The fundamental cycle of $c_1-1$ is $C'=\{c_1-1\} \cup C$ with minimal element $c_1-1$.
    Hence, $\lambda|_{C'}$ ends with an integer larger than $c_1 - 1$ which
    contradicts $\lambda\in\Lab(S^{(A)})$ by \Cref{prop:LabM-general}.
    
    (ii)
    Assume there exists a cycle $C_A=\{c_1<\cdots < c_k\}\subseteq T$ in $G_{S,A}$.
    By the construction of the graph, both edges $c_1$ and $c_k$ are incident to the root of $G_{S,A}$.
    If $C_A$ does not contain any low valley in $B\setminus A$, then $C_B=C_A$ is a cycle in $G_{S,B}$ which contradicts the assumption that $T$ is a tree in $G_{S,B}$.
    If $C_A$ contains a low valley in $B\setminus A$, let $c_j$ be the maximal one. 
    Then, since $c_j$ is incident to the root in $G_{S,B}$, the set $C_B=\{c_j< \cdots < c_k\}\subseteq T$ forms a cycle in $G_{S,B}$.
    
    (iii)
    If $C_A$ does not contain any low valley in $B\setminus A$, then $C_B=C_A$ is the fundamental cycle of $x$ in $G_{S,B}$.
    Otherwise, write $C_A=\{c_1<\cdots < c_k\}$, and let $c_j=b$ be the largest element of $C_A$ in $B\setminus A$.
    As in the proof of~(ii), the set $\{c_j< \cdots < c_k\}\subseteq T\cup\{x\}$ forms a cycle in $G_{S,B}$, since $c_j$ and $c_k$ are incident to the root in $G_{S,B}$, see \Cref{fig:containment-of-cycles}.
    Since $T$ is a spanning tree in $G_{S,B}$, this cycle is not contained in $T$.
    Hence, it contains $x$ and is the fundamental cycle $C_B$ of $x$ in $G_{S,B}$.
\end{proof}

\begin{lemma}\label{lemma:S-sequence-boolean-structure}
Let $S$ be a connected snake contained in a lattice path matroid $M$, and let $A,B\subseteq \lowvalley_M(S)$. Then
    \begin{enumerate}[label=(\roman*)]
        \item $\Lab(S^{(B)}) \subseteq \Lab(S^{(A)})$ if $A\subseteq B$, and
        \item $\Lab(S^{(A\cup B)}) = \Lab(S^{(A)}) \cap \Lab(S^{(B)}) $.
    \end{enumerate}
\end{lemma}
\begin{proof}
    (i)
    Let $\lambda=(\lambda_1,\dots , \lambda_r)\in \Lab(S^{(B)})$.
    Hence, $T_\lambda=\{\lambda_1,\dots , \lambda_r\}$ is a spanning tree of $G_{S,B}$ and, by \Cref{lemma:containment-of-cycles}~(ii), also of $G_{S,A}$.
    For an edge $x\notin T_\lambda$ let $C_A$ (resp. $C_B$) be the fundamental cycle of $x$ with respect to $T_\lambda$ in $G_{S,A}$ (resp. $G_{S,B}$).
    For the last element $\lambda_\bullet$ of $\lambda|_{C_B}$ we have $\lambda_\bullet<x$.
    By \Cref{lemma:containment-of-cycles}~(iii), every element in $C_A\setminus C_B$ is smaller than any element in $C_B$.
    Thus, the last element of $\lambda|_{C_A}$ is at most $\lambda_\bullet$ which is smaller than $x$.
    This gives $\lambda\in \Lab(S^{(A)})$.

    (ii)
    By~(i), we already have $\Lab(S^{(A\cup B)})\subseteq \Lab(S^{(A)})\cap \Lab(S^{(B)})$.
    For the reverse inclusion, let $\lambda\in\Lab(S^{(A)})\cap\Lab(S^{(B)})$.
    By~\Cref{lemma:containment-of-cycles}~(i), $T_\lambda$ is a spanning tree in $G_{S,A\cup B}$.
    For $x\notin T_\lambda$, denote its fundamental cycles in $G_{S,A}$, $G_{S,B}$ and $G_{S,A\cup B}$ by $C_A$, $C_B$ and $C$, respectively.
    By \Cref{lemma:containment-of-cycles}~(iii) for $A\subseteq A\cup B$, we have $C=C_A$, or $\min C\in B\setminus A$.
    Similarly, for $B\subseteq A\cup B$, we have $C=C_B$, or $\min C\in A\setminus B$.
    Then, since $A\setminus B$ and $B\setminus A$ are disjoint sets, $C=C_A$ or $C=C_B$.
    In both cases, $\lambda|_C$ ends with an integer smaller than $x$ by \Cref{prop:LabM-general}, and thus $\lambda\in\Lab(S^{(A\cup B)})$.
\end{proof}

Motivated by \Cref{lemma:S-sequence-boolean-structure},
we call a sequence $\lambda$ \Dfn{$(M,S)$-surviving}
if 
\begin{equation}\label{def:MS-surviving}
    \lambda\in \Lab(S) \text{ and }
    \lambda\notin \Lab(S^{(b)}) \text{ for all }
    b\in \lowvalley_M(S) \,.
\end{equation}
We remark that this definition is indeed equivalent to the one given in \Cref{subsec:intro-main-results-lpm} following the definition of $S^{(A)}$ and the characterization of internal faces.

\begin{corollary}\label{cor:Lab-incl-excl}
    The $(M,S)$-surviving sequences are given as the inclusion-exclusion of $S^{(A)}$-sequences
    \[
    \{\lambda\in \Lab(S)\mid \lambda \ (M,S)\text{-surviving}\}
    =
    \bigcup_{\substack{A \subseteq \lowvalley_M(S) \\ |A| \text{ even}}} \Lab(S^{(A)})
    \setminus
    \bigcup_{\substack{A \subseteq \lowvalley_M(S) \\ |A| \text{ odd}}} \Lab(S^{(A)})
    \]
    where both sets are considered as multisets and every element in the $|A|$-odd multiset is indeed an element of the $|A|$-even multiset.
    In particular,
    \[
        \Psi_{(M,S)}^\pm(\aaa,\bbb)
        = \sum_{\substack{\lambda\\(M,S)\text{-surviving}}} \operatorname{m}_{\Des(\lambda)}(\aaa,\bbb) \,.
    \]
\end{corollary}
\begin{proof}
    The first equation is an immediate consequence of \Cref{lemma:S-sequence-boolean-structure} and the principle of inclusion-exclusion.
    For a connected snake $S$ contained in $M$, we therefore get
    \begin{align*}
        \Psi_{(M,S)}^\pm(\aaa,\bbb)
        &= \sum_{A\subseteq \lowvalley_M(S)} (-1)^{|A|} 
        \sum_{\lambda\in \Lab(S^{(A)})} \operatorname{m}_{\Des(\lambda)}(\aaa,\bbb) \\
        &= \sum_{\substack{\lambda\\(M,S)\text{-surviving}}} \operatorname{m}_{\Des(\lambda)}(\aaa,\bbb) \,. \qedhere
    \end{align*}
\end{proof}

\begin{proof}[Proof of \Cref{thm:intro-ab-via-MS}]
    Let $M$ be a connected lattice path matroid.
    By \Cref{lma:ab-via-lowvalley}, the $\aaa\bbb$-index of $M$ decomposes as
    \[
        \Psi_M(\aaa,\bbb) = \sum_{\substack{S \text{ connected snake}\\ \text{contained in } M}} \Psi_{(M,S)}^\pm(\aaa,\bbb) \,.
    \]
    It follows with \Cref{cor:Lab-incl-excl} that
    \[
        \Psi_M(\aaa,\bbb)
        = \sum_{\substack{S \text{ connected snake}\\ \text{contained in } M}} \, \sum_{\substack{\lambda\\(M,S)\text{-surviving}}} \operatorname{m}_{\Des(\lambda)}(\aaa,\bbb) \,. \qedhere
    \]
\end{proof}

\begin{remark}[Why low valley vs. high peak matters]\label{rem:high-peak-not-boolean}
    \Cref{lemma:S-sequence-boolean-structure} does not apply to high peaks,
    since $\Lab(S_{(A)})\not\subseteq\Lab(S)$ in general.
    For example, let $M=\LPM{NNEE}{EENN}$
    and~$S=\LPM{NNEE}{ENEN}$ with high peak $2$.
    Then,
    \[
        \Lab(S)=\{(1,2),(2,1),(3,1)\}
        \quad\text{and}\quad
        \Lab(S_{(2)})=\{(1,3),(3,1)\}\,,
    \]
    and the sequence $(1,3)\in \Lab(S_{(2)})$ is not contained in $\Lab(S)$.
    \Cref{lemma:S-sequence-boolean-structure} holds for high peaks only when we consider the reversed order on the ground set, $\{n<n-1<\dots < 1\}$.
    Let $\overleftarrow{S}_{\!(A)}$ be the matroid $S_{\! A}$ with both the word for the upper path and the word for the lower path reversed.
    The matroid isomorphism $i\mapsto n+1-i$ then sends $S^{(A)}$ to a matroid isomorphic to $\overleftarrow{S}_{\!(A)}$, identifying low valleys with high peaks.
\end{remark}

\subsection{\texorpdfstring{$(M,S)$}{(M,S)}-surviving sequences}\label{subsec:MS-surviving-seq}
Let $S$ be a connected snake contained in a lattice path matroid $M$.
\Cref{thm:intro-LPM-surviving-seq} provides a combinatorial interpretation of the $(M,S)$-surviving sequences using the graph $G_S^M$ only.
Its proof is based on the following lemma, characterizing the fundamental cycles that change from $G_S^M$ to $G_{S,\{b\}}^M$ for a low valley~$b\in \lowvalley_M(S)$.

\begin{lemma}\label{lemma:lowvalley-in-fundamental-cycle}
    Denote the edges in the initial path of $G_S$ by $1=y_1<y_2<\dots<y_r$.
    Let $\lambda\in\Lab(S)$ be an $S$-sequence, and let $b=y_i\in\lowvalley_M(S)$ be a low valley.
    An edge $x=y_\ell\notin T_\lambda$ satisfies $b\in C(x,T_\lambda)$ if and only if the set $\{ y_j \geq b \mid y_j\notin T_\lambda \}$ is non-empty and $x$ is its minimum.
\end{lemma}
\begin{proof}
    ($\Rightarrow$)
    Let $x\notin T_\lambda$ be not incident to the root with $b\in C(x,T_\lambda)$.
    By \Cref{lemma:containment-of-cycles}~(i), $T_\lambda$ is a spanning tree in $G_{S,\{b\}}$.
    The fundamental cycle of $x$ in $G_{S,\{b\}}$, denoted by $C_x^{(b)}$, contains all elements in $C(x,T_\lambda)$ that are not smaller than $b$, see \Cref{lemma:containment-of-cycles}~(iii).
    Since $x\in C_x^{(b)}$ and since both $x$ and $b$ are edges in the initial path of $G_S$,
    we obtain
    \[
        \{ b=y_i < y_{i+1} < \dots < y_{i+m} =x \} \subseteq C^{(b)}_x \subseteq C(x,T_{\lambda}) \,.
    \]
    In particular, $x$ is the smallest edge in $\{ y_j \geq b \mid y_j\notin T_\lambda \}$.

    ($\Leftarrow$)
    If $b\notin T_\lambda$, then $x=b$ and $b\in C(x,T_\lambda)$.
    Otherwise, if $b\in T_\lambda$, let $x=y_\ell$ be the smallest edge not incident to the root with $x\notin T_\lambda$ and $x>b$.
    Since $T_\lambda$ contains no cycles in $G_{S,\{b\}}$, see \Cref{lemma:containment-of-cycles}~(i), every edge $z$ incident to the root and with $b<z<x$ cannot be contained in $T_\lambda$.
    Therefore, the smallest edge~$c$ in $C(x,T_\lambda)$, which is incident to the root, must be smaller than $b$, $c<b$.
    Let $d>c$ be the unique other edge in $C(x,T_\lambda)$ that is incident to the root, then $C(x,T_\lambda)$ consists of $c$, $d$, and all edges of the initial path in between:
    \[
        \tikzset{near top/.style={graph label/.append style={pos=0.25}}}
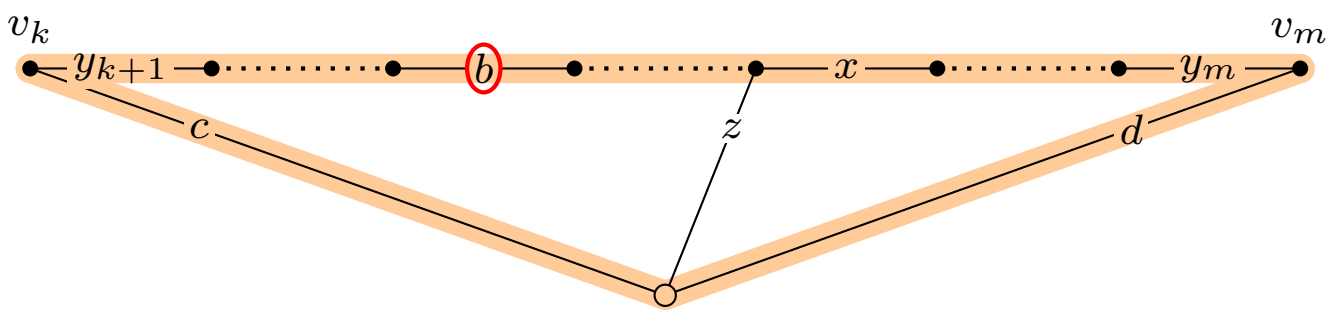
    \begin{tikzpicture}[x=1cm, y=1cm]
        \drawgraph[1/v_k, 8/v_m]{8}{
            c/1/0/0/{near top,C}, y_{k+1}/1/2/0/C, {}/2/3/0/{dots,C},
            {b}/3/4/0/{lv,C},
            z/5/0/0/near top, {}/4/5/0/{dots,C},
            {x}/5/6/0/C, {}/6/7/0/{dots,C}, y_m/7/8/0/C,
            d/8/0/0/{near top,C}}
    \end{tikzpicture}
    \]
    Since $b$ and $x$ are edges of the initial path with $c<b<x$ and $x\in C(x,T_\lambda)$, we obtain $b\in C(x,T_\lambda)$.
\end{proof}

\begin{corollary}\label{cor:Lab-Sb-condition}
    Let $\lambda\in \Lab(S)$ and let $b\in \lowvalley_M(S)$.
    Then, $\lambda\in\Lab(S^{(b)})$ if and only if
    $\lambda|_{C(x,T_\lambda)_{\geq b} }$ ends with an integer smaller than $x$ for every $x\notin T_\lambda$ not incident to the root with $b\in C(x,T_\lambda)$.
\end{corollary}
\begin{proof}
    By \Cref{lemma:containment-of-cycles}~(i), $T_\lambda$ is a spanning tree in $G_{S,\{b\}}^M$.
    Let $x\notin T_\lambda$ with fundamental cycle $C=C(x,T_\lambda)$ in $G_S^M$.
    Consider the fundamental cycle $C^{(b)}$ of $x$ in $G_{S,\{b\}}^M$. That is, by \Cref{lemma:containment-of-cycles}~(iii),
    \begin{equation*}
            C^{(b)} =
            \begin{cases}
                C & \text{if } b\notin C\,, \\
                C_{\geq b}=\{z\in C \mid z\geq b\} & \text{if } b\in C\,.
            \end{cases}
    \end{equation*}
    If $x$ is incident to the root, it must be the maximal element of $C$ since $\lambda|_{C}$ ends with an integer smaller than $x$.
    Then all other elements of $C^{(b)}$ are smaller than $x$, so the condition for $x$ holds automatically.
    If $x$ is not incident to the root and $b\in C$, then $b\leq x$ by \Cref{lemma:lowvalley-in-fundamental-cycle}.
    Together with \Cref{prop:LabM-general}, we obtain that a sequence $\lambda\in \Lab(S)$ is an $S^{(b)}$-sequence if and only if $\lambda|_{C^{(b)}}$ ends with an integer smaller than $x$ for every $x\notin T_\lambda$.
\end{proof}

\Cref{thm:intro-LPM-surviving-seq} now follows from \Cref{cor:Lab-Sb-condition}.
\begin{proof}[Proof of \Cref{thm:intro-LPM-surviving-seq}]
    An $S$-sequence $\lambda$ is $(M,S)$-surviving if and only if $\lambda\notin\Lab(S^{(b)})$ for all $b\in\lowvalley_M(S)$, hence, by \Cref{cor:Lab-Sb-condition}, if and only if, for every $b\in\lowvalley_M(S)$,
    \begin{equation*}
    \begin{array}{c}
            \text{there exists an edge $x\notin T_\lambda$ not incident to the root with $b\in C(x,T_\lambda)$}
            \\[0.5ex]
            \text{such that } \lambda|_{C(x,T_\lambda)_{\geq b}} \text{ ends with an integer larger than } x \,.
    \end{array}
    \end{equation*}
    Note that if $b\notin T_\lambda$, then $x=b$ satisfies this latter property automatically, which concludes the proof.
\end{proof}

In \Cref{subsec:MS-surviving-trees}, we determine the set of trees $T_\lambda$ for $(M,S)$-surviving sequences $\lambda$, which simplifies \Cref{thm:intro-LPM-surviving-seq} further. 

\subsection{\texorpdfstring{$(M,S)$}{(M,S)}-surviving trees}\label{subsec:MS-surviving-trees}
A spanning tree $T=\{t_1< t_2 < \dots < t_r\}$ can be written as $T=T_\lambda$ for $\lambda\in\Lab(S)$ if and only if
$(t_r>t_{r-1}>\dots> t_1)$ is an $S$-sequence, or equivalently, if and only if $T$ is a tree in $G_S$ and every fundamental cycle $C(x,T)$ has minimum smaller than $x$.
This recovers a classical correspondence between NBC-bases $T$ and decreasing maximal chains under the minimal labeling, see \Cref{rem:NBC}.
We call $T$ an \Dfn{$(M,S)$-surviving tree} if $T=T_\lambda$ for an $(M,S)$-surviving sequence $\lambda$.
Not every spanning tree $T$ is an NBC-basis nor an $(M,S)$-surviving tree.
We characterize NBC-bases of $S$ in \Cref{lemma:S-tree-iff-tree-partition},
and $(M,S)$-surviving trees in \Cref{prop:surviving-tree-characterization}.

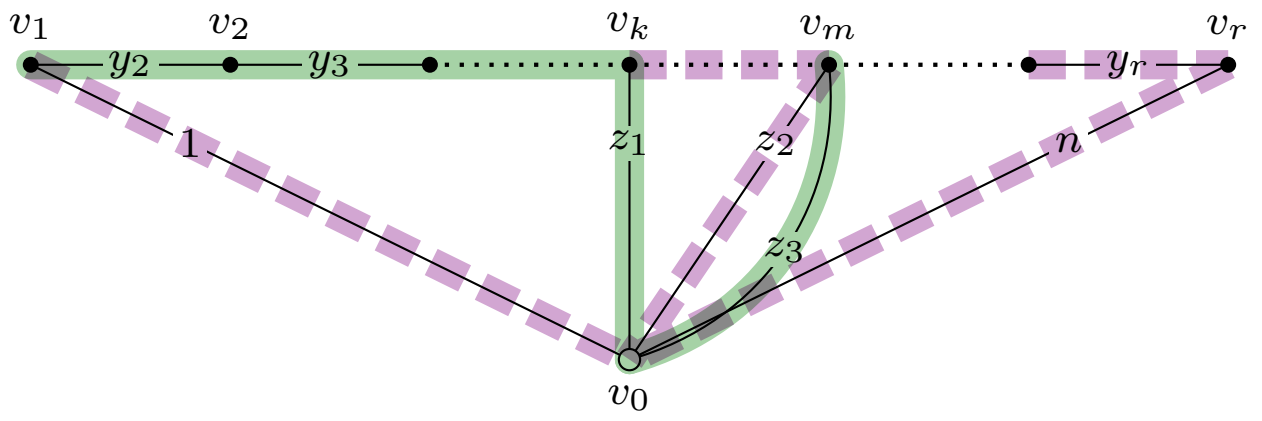
\begin{figure}[H]
    \centering
    \tikzset{near top/.style={graph label/.append style={pos=0.25}}}
    \begin{tikzpicture}[x=1.1cm, y=1.3cm]
        \drawgraph[0/v_0, 1/v_1, 2/v_2, 4/v_k, 5/v_m, 7/v_r]{7}{
            1/1/0/0/{near top,P},
            y_2/1/2/0/T, y_3/2/3/0/T, {}/3/4/0/{dots,T}, z_1/4/0/0/{near top,T},
            {}/4/5/0/{dots,P}, z_2/5/0/0/{near top,P},
            z_3/5/0/40/T,
            {}/5/6/0/dots, y_r/6/7/0/P, n/7/0/0/{near top,P}}
    \end{tikzpicture}
    \caption{The trees $\{1\}$, $T^1$, $T^2$, $T^3=\{z_3\}$, \dots, $T^{n-r}$ from~\eqref{def:tree-partition} are shaded alternately in violet (dashed) and green.}
    \label{fig:tree-partition-general}
\end{figure}

Denote the edges of $G_S$ incident to the root by $1=z_0<z_1<\dots<z_{n-r}=n$.
The set~$\{1\}$ together with the sets
\begin{equation}\label{def:tree-partition}
    T^i= \{y\in [n] \mid z_{i-1} < y \leq z_{i}\}
    \quad \text{for } 1\leq i \leq n-r 
\end{equation}
form a partition of the edges of $G_S$, see \Cref{fig:tree-partition-general}.

\begin{lemma}\label{lemma:S-tree-iff-tree-partition}
    Let $T$ be a spanning tree in $G_S$, and
    let $T^1, \dots , T^{n-r}$ be the trees in $G_S$ defined in \eqref{def:tree-partition}.
    Then, $T=T_\lambda$ for an $S$-sequence $\lambda$
    if and only if
    \[
        |T^i\setminus T|=1 \text{ for all } 1\leq i \leq n-r\,.
    \]
    In particular, the initial path of $G_S$ is such a tree.
\end{lemma}
\begin{proof}
    Let $1= z_0 < z_1 < \dots < z_{n-r}=n $ be the edges incident to the root in $G_{S}$.
    Recall that a tree $T$ equals $T_\lambda$ for some $\lambda\in \Lab(S)$ if and only if every edge $x\notin T$ is not the minimum of its fundamental cycle $C(x,T)$. 

    ($\Leftarrow$)
    Let $T$ be a spanning tree with $T^i\setminus T=\{x_i\}$ for $1\leq i\leq n-r$.
    The fundamental cycles are therefore the cycles $C(x_i,T)$ whose minimal element is smaller than $x_i$.
    Hence, $T=T_\lambda$ for some $\lambda\in \Lab(S)$.

    ($\Rightarrow$)
    Let $\lambda\in \Lab(S)$ and $1\leq i \leq n-r$.
    We have $|T^i\setminus T_\lambda|\geq 1$, since otherwise $\{z_{i-1}\}\cup T^i$ is a cycle in $G_S$ or $z_{i-1}\notin T_\lambda$ is the minimum of its fundamental cycle.
    If $z_i\notin T_\lambda$ for some $i\geq 1$, then $z_i$ is the maximum of $C(z_i, T_\lambda)$, so $T^i\subseteq C(z_i, T_\lambda)$.
    Consequently, 
    \begin{equation*}
        T^i\setminus T_\lambda=\{z_i\}
        \quad \text{whenever} \quad
        z_i\notin T_\lambda \,. \tag{$\star$}
    \end{equation*}
    For any $x\in T^i\setminus T_\lambda$ that is not incident to the root, the fundamental cycle $C(x,T_\lambda)$ must contain $T^i\setminus\{z_i\}$, that is, all edges not incident to the root between the edges $z_{i-1}$ and~$z_i$.
    By ($\star$),~$z_i$ must also be contained, which gives $T^i\setminus T_\lambda=\{x\}$.
\end{proof}

\begin{figure}[H]
    \centering
    \begin{tikzpicture}[x=1.3cm]
        \drawgraph[0/v_0, 2/v_{i-1}, 3/v_i]{5}{
            \tilde z_0/0/1/0/{T}, \tilde y_1/1/2/0/{T}, {b-1}/2/0/0/e, {b=\tilde y_2}/2/3/0/{lv,T},
            {x}/3/4/0/e, {\tilde z_1}/4/0/0/{T}, \tilde y_3/4/5/0/T, n/5/0/0/e}
    \end{tikzpicture}
    \caption{In the notation of the proof of \Cref{prop:surviving-tree-characterization}, a graph $G_S$ with tree $T=\{\tilde z_0<\tilde z_1\}\cup\{\tilde y_1<\tilde y_2<\tilde y_3\}$ (green) and low valley $b=\tilde y_2$.}
    \label{fig:proof-surviving-tree}
\end{figure}
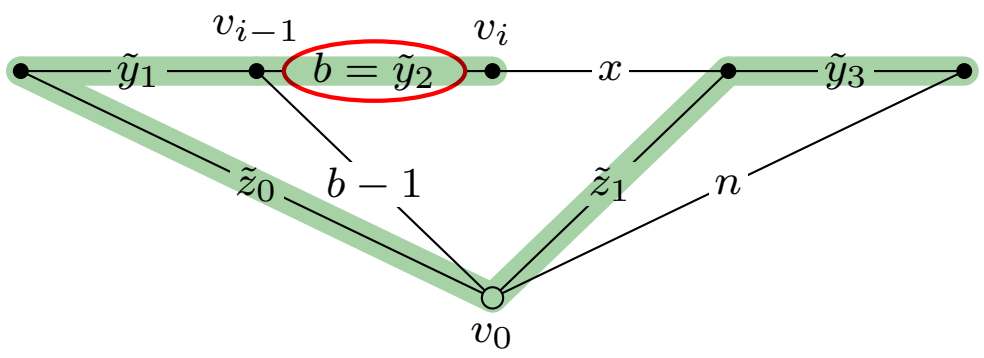

\begin{proposition}\label{prop:surviving-tree-characterization}
    A spanning tree $T$ in $G_S$ is $(M,S)$-surviving if and only if 
    \begin{enumerate}
        \item $|T^i\setminus T|=1$ for all $1\leq i \leq n-r$, and
        \item for every $b\in \lowvalley_M(S)$, there exists an edge $z\in T$ incident to the root with $z>b$.
    \end{enumerate}
\end{proposition}
\begin{proof}
    ($\Rightarrow$)
    Let $T=T_\lambda$ for an $(M,S)$-surviving sequence $\lambda$.
    The first condition follows by \Cref{lemma:S-tree-iff-tree-partition} since $\lambda\in \Lab(S)$.
    For $b\in \lowvalley_M(S)$, there exists a fundamental cycle $C(x,T_\lambda)$ for an edge $x\geq b$ not incident to the root 
    by \Cref{thm:intro-LPM-surviving-seq}, since $x=b$ works trivially when $b\notin T$.
    The maximum of this cycle $C(x,T_\lambda)$ is incident to the root, is larger than $b$ and lies in $T_\lambda$.

    ($\Leftarrow$)
    Let~$\tilde z_0 < \tilde z_1 < \dots < \tilde z_k$ be the edges in $T$ that are incident to the root, and let~$\tilde y_1 < \dots < \tilde y_\ell$ be the edges not incident to the root, see \Cref{fig:proof-surviving-tree}.
    We claim that 
    \begin{equation*}
        \lambda = ( \tilde y_1, \tilde y_{2}, \dots , \tilde y_\ell,
        \tilde z_k, \tilde z_{k-1}, \dots , \tilde z_0 ) \quad \text{is $(M,S)$-surviving.}
    \end{equation*}
    An edge $z\notin T$ incident to the root is either the minimum or the maximum of its fundamental cycle $C(z,T)$ in $G_S$.
    Condition (1) provides with \Cref{lemma:S-tree-iff-tree-partition} that $T$ is an NBC-basis, hence $z$ must be the maximal element of $C(z,T)$.
    In particular, $\lambda|_{C(z,T)}$ ends with the minimum of $C(z,T)$ which is smaller than $z$.
    
    For an edge $x\notin T$ not incident to the root, the unique two edges in $C(x,T)$ that are incident to the root are $\tilde z_{i-1}< \tilde z_i$ for some $i\in \{1,\dots , k\}$.
    In particular, $\lambda|_{C(x,T)}$ ends with $\tilde z_{i-1} < x$.
    By \Cref{prop:LabM-general}, we can conclude $\lambda\in \Lab(S)$.

    If there exists no low valleys in $(M,S)$, $\lambda$ is also an $(M,S)$-surviving sequence.
    Otherwise, let $b\in \lowvalley_M(S)$. Since $T$ is an NBC-basis, we have $1=\tilde z_0 \in T$. Together with condition (2), there exists an index $i$ such that $\tilde z_{i-1} < b < \tilde z_i$.
    Moreover, since $T$ is a spanning tree, there must exist a unique edge $x\notin T$ not incident to the root with $\tilde z_{i-1} < x < \tilde z_i$.
    The fundamental cycle $C(x,T)$ contains $\tilde z_{i-1},\tilde z_i$ and all edges between these two that are not incident to the root.
    In particular, $b\in C(x,T)$, and hence, $\lambda_{C(x,T)_{\geq b}}$ ends with $\tilde z_{i} > x$.
    It follows with \Cref{thm:intro-LPM-surviving-seq} that $\lambda$ is $(M,S)$-surviving.
\end{proof}

\begin{example}
    The graph $G_{S_4}$ from the running example is partitioned into the trees $\{1\}$, $T^1=\{2\}$, $T^2=\{3,4\}$ and $T^3=\{5,6\}$.
    Using \Cref{lemma:S-tree-iff-tree-partition}, a tree that admits an $S_4$-sequence avoids $2$, exactly one of $3,4$, and exactly one of $5,6$.
    These trees are highlighted in $G_{S_4}^M$ below.
    
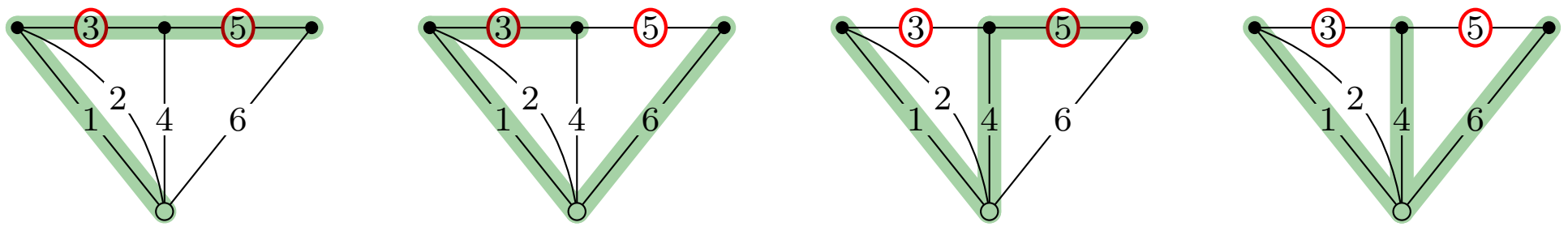
\[
        \begin{tikzpicture}
            \drawgraph{3}{1/0/1/0/T, 3/1/2/0/lvT, 5/2/3/0/lvT, 6/3/0/0/e, 2/1/0/30/e, 4/2/0/0/e}
        \end{tikzpicture}
        \qquad
        \begin{tikzpicture}
            \drawgraph{3}{1/0/1/0/T, 3/1/2/0/lvT, 5/2/3/0/lv, 6/3/0/0/T, 2/1/0/30/e, 4/2/0/0/e}
        \end{tikzpicture}
        \qquad
        \begin{tikzpicture}
            \drawgraph{3}{1/0/1/0/T, 3/1/2/0/lv, 5/2/3/0/lvT, 6/3/0/0/e, 2/1/0/30/e, 4/2/0/0/T}
        \end{tikzpicture}
        \qquad
        \begin{tikzpicture}
            \drawgraph{3}{1/0/1/0/T, 3/1/2/0/lv, 5/2/3/0/lv, 6/3/0/0/T, 2/1/0/30/e, 4/2/0/0/T}
        \end{tikzpicture}
    \]
    Among these, only the trees $\{1,\textcolor{red}{3},6\}$ and $\{1,4,6\}$ are $(M,S_4)$-surviving trees, and they provide the surviving sequences $(3,6,1)$ and $(6,4,1)$ as computed in \Cref{ex:running-intro}.
\end{example}

\section{\texorpdfstring{$\aaa\bbb$}{ab}-indices for Schubert matroids}\label{sec:applications}
A \Dfn{Schubert matroid} on $[n]$ of rank~$r$ is a lattice path matroid $M=M[U,L]$ with trivial upper path $U=\lp{N}^r\lp{E}^{n-r}$ and~$M$ is the \Dfn{uniform matroid} if also the lower path is trivial, $L=\lp{E}^{n-r}\lp{N}^r$.
In this section, we apply our results to Schubert matroids and conclude with the example of the uniform matroid.

\begin{figure}[H]
    \centering
    \tikzset{
        composition brace/.style={decorate, decoration={brace, amplitude=3pt}, thin},
        composition label/.style={font=\scriptsize, inner sep=2pt}
    }
    \begin{tikzpicture}[scale=0.65, baseline=(current bounding box.center)]
        \drawcells[snakestyle]{0/0,1/0,2/0,2/1,2/2,3/2,4/2,4/3,5/3}
        \drawpath{E,E,E,N,N,E,E,N,E,N}
        \drawpath{N,E,E,N,N,E,E,N,E,E}
        \draw[red, very thick]
            (2,1) -- (2,2)
            (4,3) -- (4,4);
        \foreach \x in {1,...,6}
            \node[font=\scriptsize, inner sep=1pt] (col\x) at (\x-0.5,4.45) {$\x$};
        \foreach \x in {2,4,6}
            \draw[red, thick] (col\x.south west) -- (col\x.north east);
    \end{tikzpicture}
    \qquad\qquad
    \begin{tikzpicture}[scale=0.65, baseline=(current bounding box.center)]
        \drawcells[snakestyle]{0/0,1/0,1/1,1/2,2/2,2/3}
        \drawpath{E,E,N,N,E,N,N}
        \drawpath{N,E,N,N,E,N,E}
        \draw[red, very thick]
            (1,1) -- (1,2)
            (2,3) -- (2,4);
        \foreach \x/\c in {1/1,2/3,3/5}
            \node[font=\scriptsize, inner sep=1pt] at (\x-0.5,4.45) {$\c$};
        \node[composition label, left] at (-0.15,0.5) {$1$};
        \draw[composition brace] (0.82,1) -- (0.82,3)
            node[midway, left=3pt, composition label] {$c_1=2$};
        \draw[composition brace] (1.82,3) -- (1.82,4)
            node[midway, left=3pt, composition label] {$c_2=1$};
    \end{tikzpicture}
    \caption{A connected snake~$S$ (left) and its compressed snake $\overline{S}=S'_c$ for $c=(2,1)\vDash 3$ (right), obtained by deleting the crossed out columns $2$, $4$ and~$6$.
    The valleys of the upper paths are marked in red.}
    \label{fig:schubert-compression}
\end{figure}
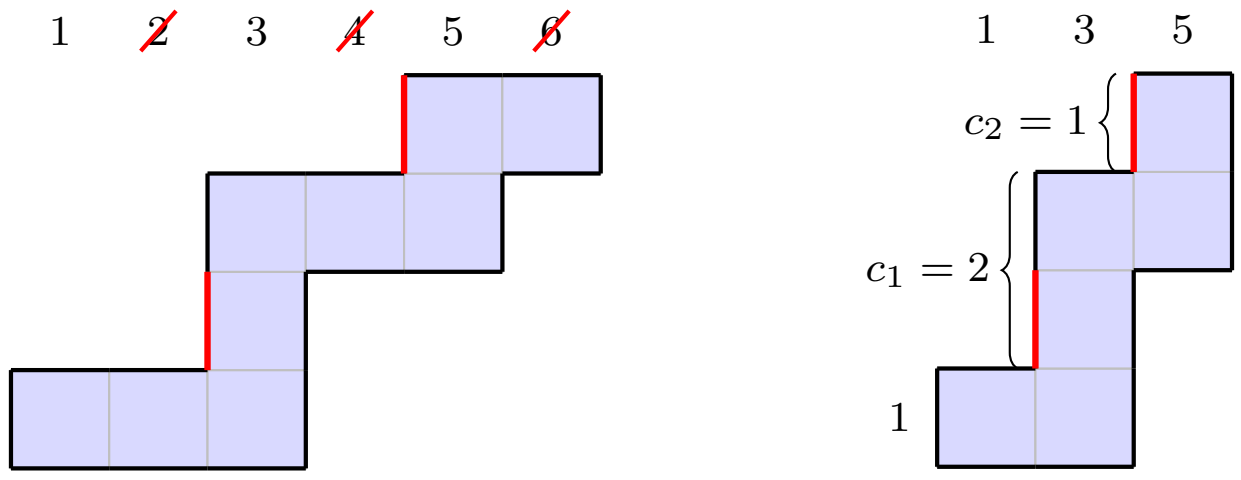

\smallskip
For a connected snake $S$, delete every column in the diagrams of $S$ that is not the first column, and that contains only a single box of $S$, to obtain the compressed snake $\overline S$, see \Cref{fig:schubert-compression}.
This compressed snake $\overline{S}$ either equals $S_c$ or $S_c'$ for a unique composition $c\vDash r-1$, see \Cref{fig:snakes-Sc}.

\begin{lemma}\label{lma:schubert-ab-via-composition}
    Let $S$ and $S'$ be two connected snakes contained in a Schubert matroid $M$ such that~$\overline{S} = \overline{S'}$.
    Then, $\Psi_{(M,S)}^{\pm}(\aaa,\bbb)=\Psi_{(M,S')}^{\pm}(\aaa,\bbb)$.
\end{lemma}
\begin{proof}
    Let $\overline{M}$ be the unique uniform matroid for which $\overline{S}$ is contained in.
    We prove that~$\Psi_{(M,S)}^{\pm}(\aaa,\bbb)=\Psi_{(\overline{M},\overline{S})}^{\pm}(\aaa,\bbb)$, from which the lemma follows.
    
    Let $S$ have a single box present in column $i>1$ of its diagram.
    Then, the edges $i-1,i$ in $G_S$ are incident to the root,
    and any $S$-sequence must not contain $i$.
    Deleting column $i$ entirely yields a connected snake whose graph is obtained from $G_S$ by deleting edge $i$ and shifting all edge labels larger than $i$ by $-1$.
    Doing this process iteratively yields $\overline S$ whose NBC-bases are obtained from those of $S$ and an order-preserving map $[n]\to[n]$.
    Since the deletion of columns only shifts the labels of low valleys in $(M,S)$, we conclude the proof. 
\end{proof}

\begin{figure}[H]
    \centering\begin{tabular}{c@{\qquad}c@{\qquad}c}
         \begin{tikzpicture}[scale=0.6]
             \drawcells[snakestyle]{0/0, 0/1, 1/1, 2/1, 2/2, 3/2}
             \drawcells{0/0, 0/1, 0/2, 1/0, 1/1, 1/2, 2/1, 2/2, 3/1, 3/2}
             \drawpath{N,N,N,E,E,E,E}
             \drawpath{E,E,N,E,E,N,N}
             \draw[red, very thick] (2,2) -- (2,3);
         \end{tikzpicture} &
         \begin{tikzpicture}[scale=0.6]
             \drawcells[snakestyle]{0/0, 0/1, 1/1, 1/2, 2/2, 3/2}
             \drawcells{0/0, 0/1, 0/2, 1/0, 1/1, 1/2, 2/1, 2/2, 3/1, 3/2}
             \drawpath{N,N,N,E,E,E,E}
             \drawpath{E,E,N,E,E,N,N}
             \draw[red, very thick] (1,2) -- (1,3);
         \end{tikzpicture} &
         \begin{tikzpicture}[scale=0.6]
             \drawcells[snakestyle]{0/0, 0/1, 1/1, 1/2}
             \drawcells{0/0, 0/1, 0/2, 1/0, 1/1, 1/2}
             \drawpath{N,N,N,E,E}
             \drawpath{E,E,N,N,N}
             \draw[red, very thick] (1,2) -- (1,3);
         \end{tikzpicture} \\
         \begin{tikzpicture}
             \drawgraph{3}{1/0/1/0/e, 2/1/2/0/e, 5/2/3/0/lv, 3/2/0/-20/e, 4/2/0/20/e, 6/3/0/0/e, 7/3/0/35/e}
         \end{tikzpicture} &
         \begin{tikzpicture}
             \drawgraph{3}{1/0/1/0/e, 2/1/2/0/e, 4/2/3/0/lv, 3/2/0/0/e, 5/3/0/-25/e, 6/3/0/10/e, 7/3/0/45/e}
         \end{tikzpicture} &
         \begin{tikzpicture}
             \drawgraph{3}{1/0/1/0/e, 2/1/2/0/e, 4/2/3/0/lv, 3/2/0/0/e, 5/3/0/0/e}
         \end{tikzpicture}
    \end{tabular}
    \caption{Two connected snakes in the same Schubert matroid (left and middle) and with the same compressed snake $S_c$ for $c=(1,1)$ (right), together with their graphs below. Low valleys are marked in red.}
    \label{fig:schubert-delete-parallel-edges}
\end{figure}
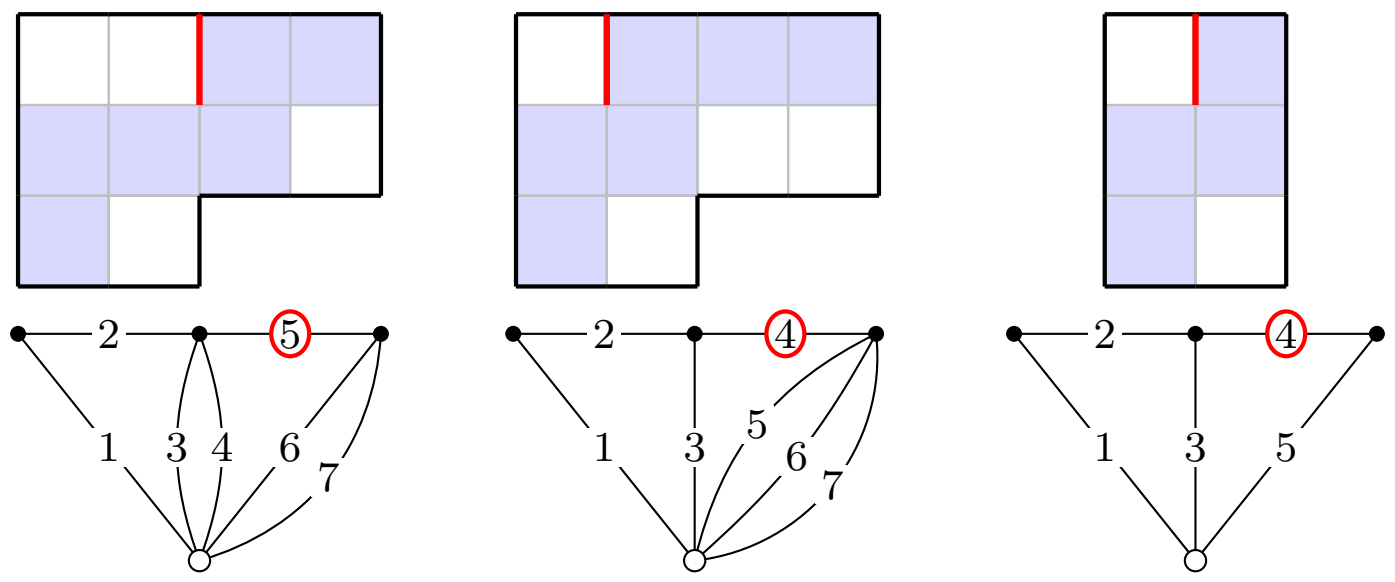

Since every valley of the upper path of a connected snake contained in a Schubert matroid is a low valley
and by~\Cref{lma:schubert-ab-via-composition},
$\Psi_{(M,S)}^\pm(\aaa,\bbb)$ only depends on the compressed snake $\overline{S}$.
For a composition $c=(c_1,\dots , c_\ell)$ of length $\ell(c)=\ell$, 
we thus define
\[
    \Psi_{c}^{\pm}(\aaa,\bbb) = \Psi_{(M,S_c)}^{\pm}(\aaa,\bbb)
        \qquad \text{and} \qquad
    \Psi_{1,c}^{\pm}(\aaa,\bbb) = \Psi_{(M',S'_c)}^{\pm}(\aaa,\bbb) \,,
\]
where $M$ and $M'$ are Schubert matroids containing $S_c$ and $S_c'$, respectively.
Note that~$S_c$ has $\ell-1$ low valleys and $S_c'$ has $\ell$ low valleys.

For a composition $c$ and a lattice path $L$, let $q_j$ be the number of east steps of $L$ at height at most $\sum_{i< j} c_i$ for $j=0,\dots , \ell(c)$, let $B_j = \min\{q_k - k \mid j\leq k \leq \ell(c) \}$, and define
\[
    f_{c}(L)=\det_{2\leq i,j \leq \ell}\binom{B_i+1}{j-i+1}, \qquad
    f'_{c}(L)=\det_{1\leq i,j \leq \ell}\binom{B_i}{j-i+1},
\]
where $\binom{a}{k}=0$ for $a<k$ and the determinant of the empty matrix is $1$.

\begin{proof}[Proof of \Cref{thm:ab-schubert}]
    If $f_{c}(L)$ and $f_{c}'(L)$ are the numbers of connected snakes $S$ contained in~$M$ with $\overline{S}=S_c$ and $\overline{S}=S_c'$, respectively,
    we obtain 
    \[
        \Psi_M(\aaa,\bbb) 
        = \sum_{c\, \vDash \, r-1} \Big( 
            f_{c}(L)\, \Psi^{\pm}_{c}(\aaa,\bbb) 
            + f_{c}'(L)\, \Psi^{\pm}_{1,c}(\aaa,\bbb)\Big) 
    \]
    by \Cref{lma:ab-via-lowvalley,lma:schubert-ab-via-composition},
    and the theorem follows with \Cref{cor:Lab-incl-excl}.
    It thus remains to prove the explicit formulas of $f_{c}(L)$ and $f_{c}'(L)$ for a composition $c=(c_1,\dots , c_\ell)$ of $r-1$ and the given lower path $L$ of $M$.

    Fix a composition $c=(c_1,\dots , c_\ell)\vDash r-1$.
    A connected snake $S$ with~$\overline{S}=S_c'$ is then determined by the positions 
    \[
        2 \leq s_1<\dots<s_\ell\leq n-r
        \qquad \text{where the column $s_j$ contains $c_j+1$ boxes.}
        \tag{$\dagger$}
    \]
    We number the rows of the diagram of~$M$ from bottom to top by $1,\dots,r$, so that row~$k$ lies between the heights $k-1$ and~$k$.
    The $c_j+1$ boxes in column $s_j$ are in the rows $p_{j-1}+1,\dots,p_j+1$ and every other column contains a single box:
    
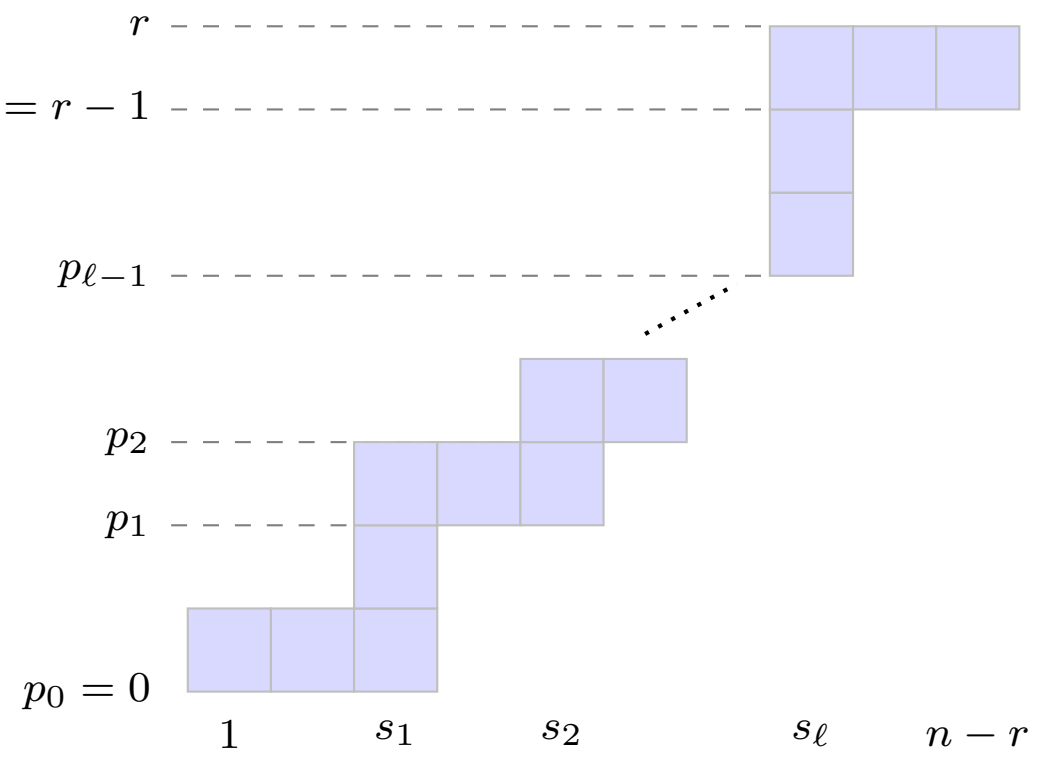
\[
        \begin{tikzpicture}[scale=0.55, every node/.style={font=\scriptsize}]
            \foreach \y/\t in {0/{p_0=0},2/{p_1},3/{p_2},5/{p_{\ell-1}},7/{p_\ell=r-1},8/{r}}
                \node[left] at (-0.2,\y) {$\t$};
            \draw[gray, dashed] (-0.2,2) -- (3,2);
            \draw[gray, dashed] (-0.2,3) -- (5,3);
            \draw[gray, dashed] (-0.2,5) -- (7,5);
            \draw[gray, dashed] (-0.2,7) -- (8,7);
            \draw[gray, dashed] (-0.2,8) -- (7,8);
            \drawcells[snakestyle]{0/0,1/0,2/0,2/1,2/2,3/2,4/2,4/3,5/3}
            \drawcells[snakestyle]{7/5,7/6,7/7,8/7,9/7}
            \draw[dotted, thick] (5.5,4.3) -- (6.6,4.9);
            \foreach \x/\t in {0/{1},2/{s_1},4/{s_2},7/{s_\ell},9/{n-r}}
                \node[below] at (\x+0.5,-0.1) {$\t$};
        \end{tikzpicture}
    \]
    The lower path $L$ of the Schubert matroid $M$ admits a unique east step in every column, and the snake~$S$ is contained in~$M$
    if and only if, in every column, the boxes of~$S$ do not lie below the east step of~$L$ in this column.
    It moreover suffices to verify this condition for columns $s_j$,
    since the heights of the east steps of~$L$ weakly increase and are at most $r-1$.
    The columns $s_{j-1}+1,\dots,s_j$ start at height~$p_{j-1}$, where $s_0=0$, and the columns after~$s_\ell$ start at height $p_\ell=r-1$.
    Therefore, $S$ is contained in $M$ if and only if, for every $1\leq j\leq \ell$, the east step of~$L$ in column~$s_j$ has height at most $p_{j-1}$, or equivalently, if and only if
    \[
        s_j\leq q_j \quad\text{for all } 1\leq j\leq \ell\,.
    \]
    We change variables by defining $t_j = s_j -j$ and $B_j = \min\{q_k -k \mid j\leq k \leq \ell \}$.
    In this way,
    \[
        1 \leq t_1 \leq \dots \leq t_\ell \leq n-r-\ell
        \qquad \text{and} \qquad
        B_1\leq \dots \leq B_\ell
    \]
    are both weakly increasing and $S$ is contained in $M$ if and only if $t_j \leq B_j$ for all $1\leq j \leq \ell$.
    If $B_1<1$, there exists no such snake and $f_{c}'(L)=1$.
    Otherwise, $f_{c}'(L)$ counts lattice paths from $(0,0)$ to $(\ell,B_\ell-1)$ that stay weakly below $\lp{N}^{B_1-1}\,\lp{E}\,\lp{N}^{B_2-B_1}\,\lp{E}\cdots\lp{N}^{B_\ell-B_{\ell-1}}\,\lp{E}$, where the $j$-th east step has height $t_j-1$.
    The result follows from \cite[Theorem 10.7.1]{krattenthaler2015} which counts lattice paths between two paths.

    The proof is analogous for $f_{c}(L)$ with $s_1=1$ in~$(\dagger)$.
    We then count the lattice paths from $(0,0)$ that end at the same point as
    $\lp{E}\,\lp{N}^{B_2}\,\lp{E}\,\lp{N}^{B_3-B_2}\,\lp{E}\cdots\lp{N}^{B_\ell-B_{\ell-1}}\,\lp{E}$
    and stay weakly below it, and apply \cite[Theorem 10.7.1]{krattenthaler2015}.
\end{proof}

\Cref{thm:ab-schubert} and its proof are not needed to compute the multiplicities $f_{c}(L)$ and $f_{c}'(L)$ for uniform matroids,
but we include this computation in the following example to illustrate a simple case.

\begin{example}
    The uniform matroid $M=M[U,L]$ on $[n]$ and of rank $r$ has lower path $L=\lp{E}^{n-r}\lp{N}^r$.
    In the notation from \Cref{thm:ab-schubert} and its proof,
    we have $q_j=n-r$ and $B_j=n-r-\ell(c)$ for every $j\geq 1$,
    hence,
    $f_{c}(\lp{E}^{n-r}\lp{N}^r)$ is the number of lattice paths from $(0,0)$ to $(\ell(c)-1,n-r-\ell(c))$, and $f_{c}'(\lp{E}^{n-r}\lp{N}^r)$ is the number of lattice paths from $(0,0)$ to $(\ell(c),n-r- \ell(c)-1)$.
    That is,
    \[
        f_{c}(\lp{E}^{n-r}\lp{N}^r)=\binom{n-r-1}{\ell(c)-1} 
        \qquad \text{and} \qquad
        f_{c}'(\lp{E}^{n-r}\lp{N}^r)=\binom{n-r-1}{\ell(c)} \,. 
    \] 
\end{example}

\printbibliography
\end{document}